\documentclass[12pt,a4paper]{article}
\usepackage{t1enc}
\usepackage[latin1]{inputenc}
\usepackage[english]{babel}
\usepackage{amssymb}
\usepackage{epsfig}
\usepackage{color}
\usepackage{amsthm, amssymb}
\usepackage{epsfig,multicol}   
\usepackage{graphics}
\newtheorem{remark}[subsubsection]{Remark}

\newtheorem{diagram}[subsubsection]{Fig.}

\newtheorem{formula}[subsubsection]{}
\newcommand{\5}{\vskip 5pt}

\newcommand\h{{\cal H}}

\begin{document}
\def\hpic #1 #2 {\mbox{$\begin{array}[c]{l} \epsfig{file=#1,height=#2}
\end{array}$}}

\def\vpic #1 #2 {\mbox{$\begin{array}[c]{l} \epsfig{file=#1,width=#2}
\end{array}$}}

\title{In and around the origin of quantum groups.}
\author{Vaughan F.R. Jones
\thanks{Supported in part by NSF Grant DMS93--22675, the Marsden fund UOA520, and the Swiss
National Science Foundation.}
} \maketitle
\abstract
Quantum groups were invented largely to provide solutions of the
Yang-Baxter equation and hence solvable models in 2-dimensional
statistical mechanics and one-dimensional quantum mechanics. They
have been hugely successful. But not all Yang-Baxter solutions
fit into the framework of quantum groups. We shall explain how
other mathematical structures, especially subfactors, provide a language and examples for
solvable models. The prevalence of the Connes tensor product 
of Hilbert spaces over von Neumann algebras leads us to speculate
concerning its potential role in describing entangled or interacting
quantum systems.
\newpage
\section{The representations of $SU(2)$}
Since $SU(2)$ is compact, any continuous representation on Hilbert space
is unitarizable and the direct sum of a family of irrecuible representations,
all of which are finite dimensional. The irreducible unitary representations
(henceforth called ``irreps'')
are easy to classify. There is exactly one of each dimension $n$ which
is often written $n=2j+1$ where $j$ is the ``spin'' of the representation.
Let $V_j$ be the vector space of the spin $j$ irrep.
Explicitly, $V_j$ can be constructed from the identity representation
on $\mathbb C^2$ as the symmetric algebra of $\mathbb C^2$.
 That is to say that $SU(2)$ acts on homogeneous polynomials of two variables $x$ and $y$
of degree $2j+1$ by extending the action $x\mapsto ax+by$, $y\mapsto cx+dy$
for a matrix
$\left( \begin{array}{cc}
a & b \cr
c & d \cr
\end{array} \right)
$ 
in $SU(2)$.

\subsection{Clebsch-Gordon rules }\label{cgr}
The tensor product decomposition for the irreps of $SU(2)$ is known as
the Clebsch-Gordon rule and is simply the following:
$$V_j\otimes V_k = \oplus_{i=|j-k|}^{j+k} V_i$$
where the equation is as $SU(2)$-modules and $i$ goes in steps of $1$.
This decomposition is easy to prove. Observe that the circle sugroup
of diagonal matrices $\left( \begin{array}{cc}
e^{i\theta} & 0 \cr
0 & e^{-i\theta} \cr
\end{array} \right)
$ 
acts in $V_j$ by diagonal matrices with respect to the basis of
monomials with eigenvalues $\{z^{2j},z^{2j-2},..., z^{-2j}\}$ (where
$z=e^{i\theta}$. These eigenvalues are the ``weights'' of the
representation. It is clear then that $V_j\otimes V_k$ has highest 
weight $z^{2j+2k}$ with multiplicity one so there is exactly one 
copy of $V_{j+k}$. Orthogonal to it we see the weight $z^{2j-2}$
with multiplicity one. Continuing in this way we are done.

When $k=2$ the Clebsch-Gordon rules say that
 $V_j\otimes V_{1/2} =V_{j+1/2}\oplus V_{j-1/2}$. Since any 
irrep is contained in a tensor power of $V_{1/2}$ one may
show that this rule alone suffices to determine all the 
Clebsch-Gordon rules. We may represent this rule graphically
as follows:

\vskip 20pt

\vpic{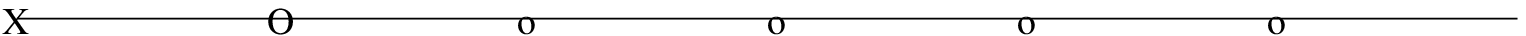} {3.5in}

\vskip 20pt
\noindent Here the vertices of the graph, known as $A_\infty$, represent 
the irreps of $SU(2)$ and an edge between two vertices means that
the irrep of one is contained in the tensor product of the 
other with $V_{1/2}$. This procedure for associating graphs with
the irreps of an object, with one privileged one, is obviously quite
general and we will use it without further explanation below. Note that
if there were multiplicity in the decomposition, one would use 
multiple edges in the graph.

\subsection{Decomposition of the tensor powers of irreps.}\label{decomposition}

If  $\pi$ is a representation of the group $G$ on the 
vector space $V$, one looks first for proper subspaces of
$V$ which are invariant under $\pi_g$ for all $g\in G$.
If $V$ is a Hilbert space and $\pi$ is unitary it is 
natural to ask that the subspace be closed, hence also 
a Hilbert space. Moreover closed subspaces of Hilbert space
are the same as projection operators-continuous linear maps
$p:\h \rightarrow \h$ with $$p=p^* \hbox{ and } p^2=p.$$
To say that the subspace is invariant is the same as saying that
the corresponding projection commutes with $\pi_g$ for all
$g\in G$. Thus the various ways in which a unitary 
representation decomposes are described entirely by projections
that commute with the group. But the set of \underline{all} 
continuous operators which commute with the group has 
the structure of an algebra to which many more techniques can
be brought to bear than on the set of its projections. 
Indeed we have just given one of the equivalent definitions
of a {\bf von Neumann algebra}, namely the algebra of 
operators commuting with a unitary group representation.
\5
If $\pi$ is any representation of any group $G$ on the vector
space $V$, there is \emph{always} a canonical algebra of linear 
transformations of $\otimes^kV$ commuting with $\otimes^k \pi$. That is the
algebra generated by the permuation group $S_k$ acting by
permuting the various tensor product components (i.e. if $\sigma$
is a permutation then $\sigma(v_1\otimes v_2 \otimes... \otimes v_k)=
v_{\sigma(1)}\otimes v_{\sigma(2)}\otimes.... \otimes v_{\sigma(k)}$
 --  or is it $\sigma^{-1}$?...).
Since the permutation group is generated by its transpositions,
this algebra is generated by $S_{12},S_{23},...S_{(k-2)(k-1)}$ where 
$S:V\otimes V\rightarrow V\otimes V$ is the map $S(v\otimes w)=w\otimes v$,
and for the rest of this paper we make the convention that if
$R:V\otimes V \rightarrow V\otimes V$ is any linear map then
for $1\leq i \leq k-1$ the linear map $R_{i(i+1)}:\otimes^kV \rightarrow \otimes^k V$
is defined by
$$R_{i(i+1)}(v_1\otimes v_2\otimes... v_i \otimes v_{i+1}...\otimes v_k)
=v_1\otimes v_2\otimes...R(v_i\otimes v_{i+1})... \otimes v_k$$

Thus one may decompose the tensor powers of $\pi$ according to 
the irreps of the symmetric group by projecting on to the subspace
of vectors (the so-called ``isotypical component'') of vectors that
transform according to that representation of $S_k$. Thus the symmetric powers
of $\pi$ are given by the trivial irrep and the antisymmetric powers
by the parity irrep. It is a well known result, sometimes called 
``Schur-Weyl duality'', that if $V=\mathbb C ^n$ and $G=SU(n)$ then
the algbebra generated by $S_k$ is in fact the algebra of all
operators commuting with $G$.

\section{The McKay correspondence}
This is a relation between closed subgroups of $SU(2)$ and
the extended Coxeter-Dynkin diagrams $\tilde A ,\tilde D$, $\tilde E$
$\tilde A_{-\infty,+\infty}$ and $\tilde D_\infty$ drawn
below (and of course $A_\infty$ drawn above).
\vskip 10pt
\vpic{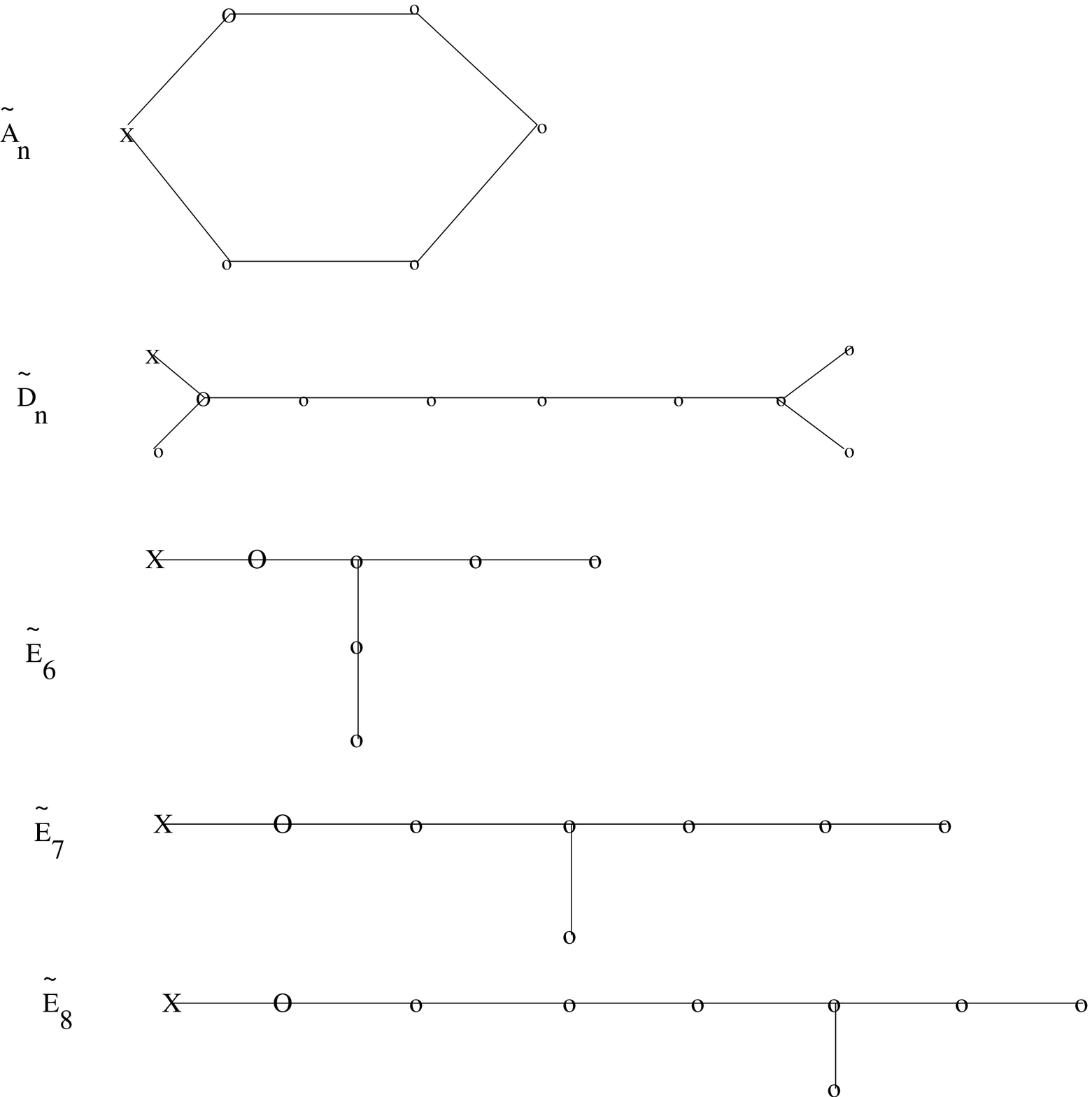} {4in}
\vskip 1pt
$A_{-\infty,\infty}$\vpic{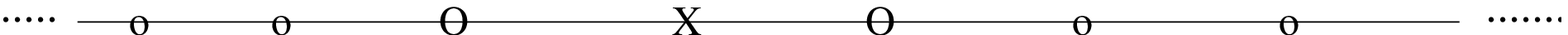} {4in}

\5 \5
$D_\infty$\vpic{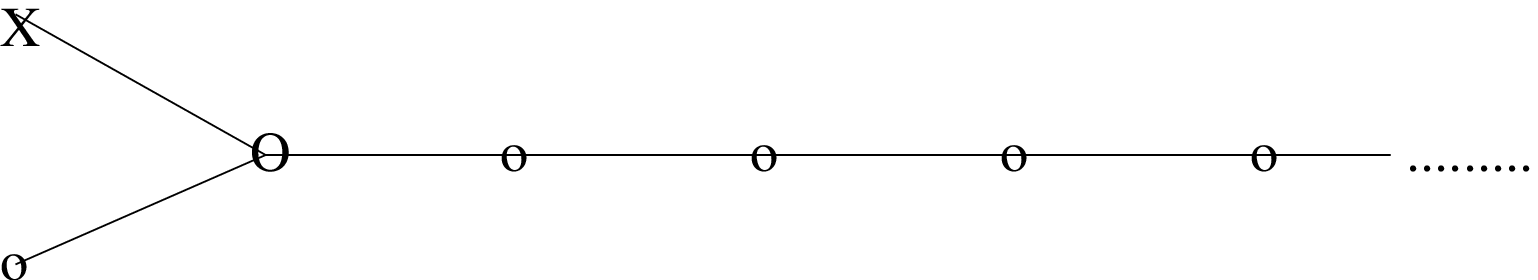} {3in}
\vskip 10pt
Let us start with the case of a finite subgroup $G$. By passing to the
quotient $SO(3)$ we see that $G$ is the double cover of either a cyclic group,
a dihedral group or the symmetry group of one of the Platonic solids-the
tetrahedron, cube/octahedron and the icosahedron/dodecahedron. 
We now form a graph for $G$ as we did for $SU(2)$ 
in \ref{cgr}. The vertices of the graph are the set of irreps of $G$ and
there are $k$ edges between two irreps if the tensor product of
one by the two dimensional
identity representation of $G$ contains $k$ copies of the other. (In fact no
multiplicity higher than one occurs here.)
The McKay correspondence asserts that the graph obtained is necessarily 
an extended Coxeter-Dynkin diagram according to the following scheme.

\begin{eqnarray*}
 \tilde{A_n}&\leftrightarrow& \hbox{Cyclic Group} \\
\tilde{D_n} &\leftrightarrow& \hbox{Dihedral Group} \\
\tilde{E_6} &\leftrightarrow&\hbox{Tetrahedral Group} \\
  \tilde{E_7} &\leftrightarrow& \hbox{Cube/Octahedron Group}\\
\tilde{E_8} &\leftrightarrow& \hbox{Icosahedral/Dodecahedral Group} 
\end{eqnarray*}

There are three infinite closed subgroups of
$SU(2)$. They are $SU(2)$ itself, the circle group $\mathbb T$ 
and the infinite dihedral group $\mathbb T \mathbb o \mathbb Z /2\mathbb Z$.
They correspond to the diagrams $A_\infty$, $A_{-\infty,\infty}$ and
$D_\infty$ respectively. Here $A_\infty$ is the graph of the Clebsch-Gordon
rules for $SU(2)$ and  $A_{-\infty,\infty}$ and $D_\infty$ are as above.
 \5
For the lovers of the empty set we must mention the trivial group 
consisiting of the identity element. It has one irreducible 
representation which, on tensoring with the identity representation
gives 2 copies of itself. So the graph of the McKay correspondence
could be taken as the graph with one vertex and two edges connecting
that vertex to itself...

The cyclic group $\mathbb Z/{n\mathbb Z}$ case requires a certain amount of care 
as the representation is not irreducible so corresponds actually
to both vertices on the graph adjacent to the trivial representation.
The cyclic groups exist as honest subgroups of $SU(2)$ and as such 
they give rise to $\tilde A_n$'s. As subgroups of $SO(3)$ they
are double covered in passing to $SU(2)$ and  
what happens depends on the parity of $n$.
We leave the somewhat confusing details as an
exercise. 
\5
The guiding light here is that the graph \underline{must}
somehow be made up from  extended $ADE$ diagrams as there is a third ingredient of 
the McKay correspondence which is to $p\times q$ matrices with
non-negative integer entries whose norm is equal to 2. (The 
norm of a matrix $\Lambda$ is the largest stretching factor for unit
vectors, or alternatively the square root of the largest eigenvalue of $\Lambda^T\Lambda$.)
In this correspondence one takes a bipartite graph with $n$ vertices, with 
disjoint vertex subsets
$X$ and $Y$ not connected to themselves, but $n=\#(X)+\#(Y)$, and constructs the 
matrix with columns labelled by $X$ and rows labelled by $Y$.
Under certain irreducibility assumptions, if the resulting matrix has norm $2$, the
graph has to be an extended $ADE$ graph.The importance of norm 2 is
explained as follows. From $A$ form the square matrix

$$\Omega =
\left( \begin{array}{cc}
0&\Lambda  \cr
\Lambda^T&0 \cr
\end{array} \right)
$$ 
(which is actually the adjacency matrix of the graph in the usual sense). The
norm of $\Omega$ is the same as that of $\Lambda$ and the Perron Frobenius
theorem on matrices with non-negative entries implies that the norm of
$\Omega$ is the eigenvalue of the unique eigenvector with positive entries.
It suffices to exhibit such a vector (whose representation theoretic nature
we will describe) for the $ADE$ diagrams to show they have norm $2$.

\5
 In the other direction one may see an a priori connection with root systems
for Lie algebras by forming $2 -\Omega$. Given that the norm of $\Omega$ is
equal to 2 and $\Omega$ is symmetric, $2-\Omega$ is positive semidefinite
so has a symmetric (real) square root $\Delta$. The relation $\Delta^2=\Omega$
says precisely that the rows of $\Delta$ are vectors which are all of length
$\sqrt{2}$ and are either orthogonal or at an angle of $120^\circ$ to
each other. Since $2$ is actually an eigenvalue of $\Omega$, the rows
of $2-\Omega$ only
span a subspace of dimension $n-1$. Up to this detail we are dealing with
a root system. In fact if any vertex of the graph is removed the resulting
set of vectors will indeed be a root system all of whose roots have the
same length. Thus we expect to see the $ADE$ Coxeter-Dynkin diagrams. 
The details are left as an exercise.

We would like to mention an amusing check on all this stuff. From the 
point of view of $SU(2)$, the reason the matrix has to have norm $2$ 
is that tensoring a representation by the identity representation 
multiplies the dimension of the representation by $2$ so that 
\emph{the vector whose entries are the dimensions of the irreps
of the closed subgroup G of $SU(2)$} is an eigenvector for
$\Omega$ of eigenvalue $2$. Conversely, if we take the Perron-Frobenius
eigenvector for $\Omega$ and normalise it so that the component corresponding
to the trivial representation is $1$, the other entries must all be 
integers, indeed they must be the dimensions of the irreps of $G$ !
We illustrate with the Perron-Frobenius eigenvector for $\tilde E_8$ below:
\5 \5
\vpic{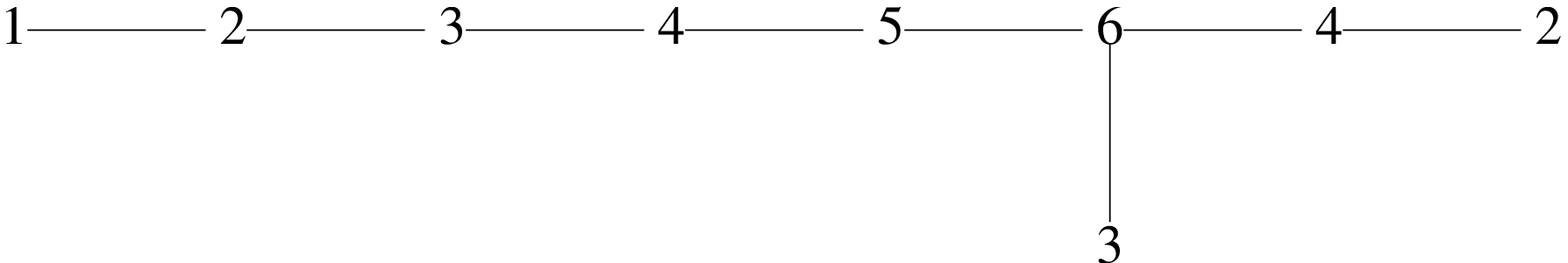} {4in}

\5
Note that the sum of the squares of the dimensions is 120. The order
of the group of rotational symmetries of the dodecahedron is 
obviously 60. The factor of 2 is due to the double covering when 
passing from $SU(2)$ to $SO(3)$.

\5
A curious question arises out of our McKay correspondence. Why did
only the \underline{extended} Coxeter-Dynkin diagrams arise? Are there naturally
arising structures whose representations are the vertices of
an ordinary $ADE$ diagram and for which the tensor product rule
can be interpreted as above? If such structures exist is there
a context in which they appear just as naturally as the McKay 
correspondence? The answer to these questions is provided by 
subfactors as we shall see.

A more obvious way to extend the McKay correspondence is to do the
same thing for $SU(3)$ and beyond. One will not of course obtain
the $ADE$ diagrams but rather graphs of norm 3, 4 and so on. Moreover
the graphs will have to be directed. The reason for undirected graphs
for $SU(2)$ and its subgroups is that the identity representation is self-conjugate.
If we had considered $U(2)$ instead we would have had to use directed 
graphs and would have found graphs with loops and directed edges.
As a very simple example for $U(3)$ here is the directed graph (of norm 3
of course) resulting from a copy of $\mathbb Z/5\mathbb Z$ in $U(3)$:
\5 \5 \5
\vpic{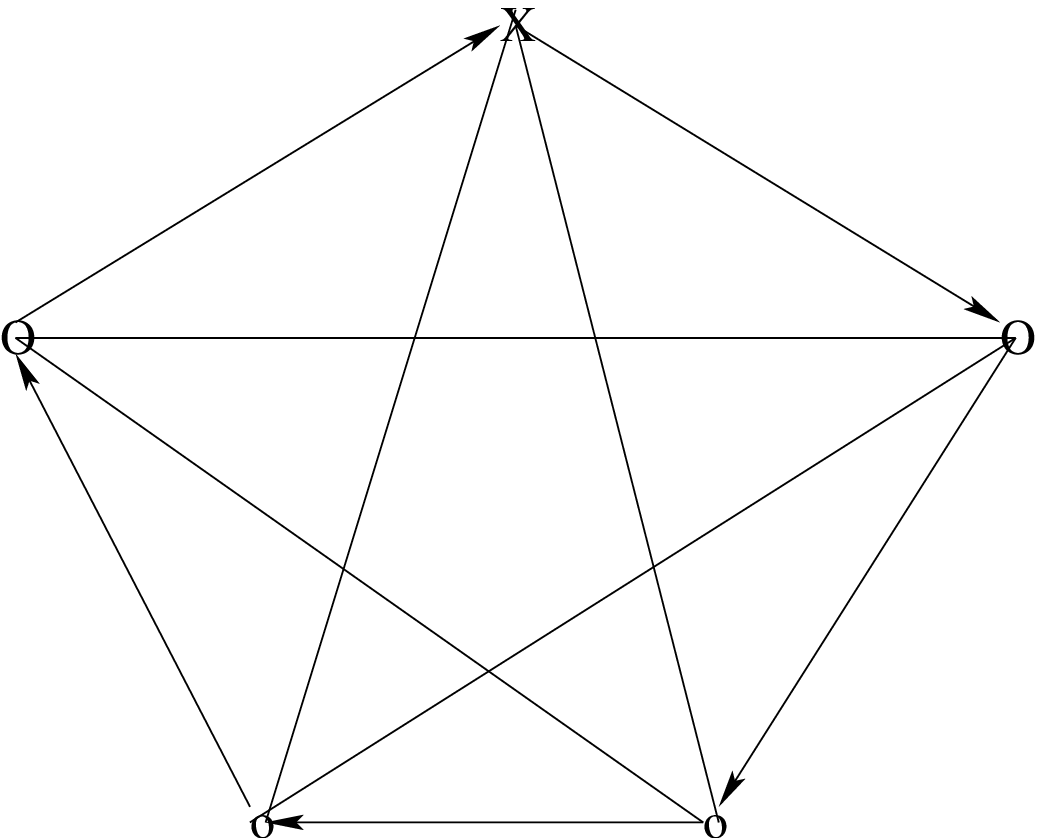} {3in}

\section{Commuting transfer matrices, the Yang-Baxter equation}
\subsection{Generalities}
In statistical mechanics systems are sometimes modelled by specifying
a set of states $\{\sigma\}$ arising from a collection of locally interacting 
sites placed on some lattice. An energy is assigned to each state
according to the model. If just a finite subset $X$ of the lattice,with $N$ sites, 
is considered, the number of states
may be finite and the ``partition function'' for $X$ is 
$$Z_X=\sum_{\sigma \in \{\hbox{states of } X\}}e^{-E(\sigma)\over
kT}$$
Some attention will have to be given to the boundary of $X$ to properly 
define $Z_X$. In general we will consider an increasing family of 
subsets $X$ whose union is the whole system.

For instance the simplest of all such models is the Ising model where the
lattice is $\mathbb Z ^n$ and $X$ is a product of intervals, depicted below
for instance when $n=2$ and $X$ is a $6\times 6$ square:
\5
\hspace{1in}\vpic{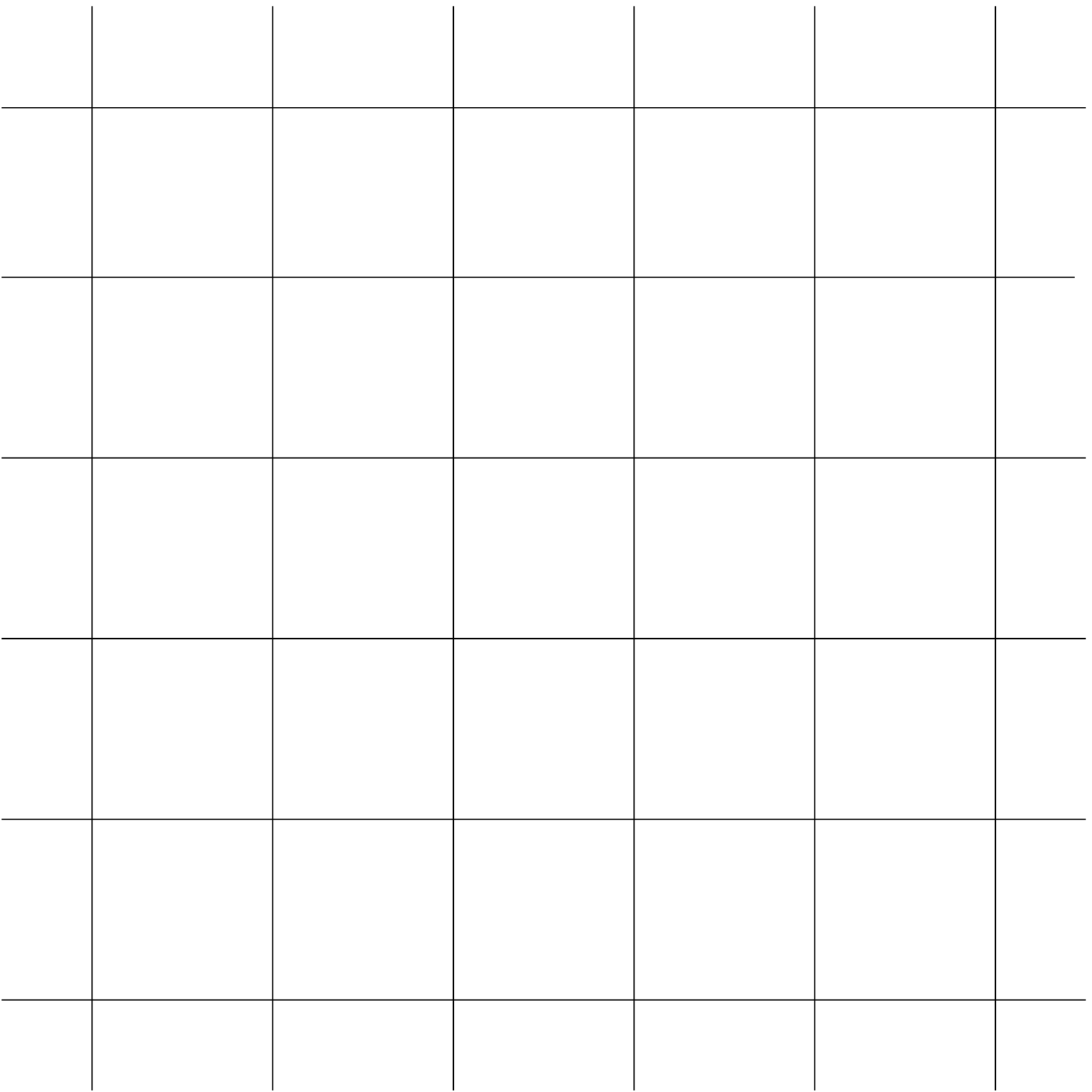} {2in}
\5

A state of the system is specified by assigning one of two ``spin'' states
$\uparrow $ and $\downarrow $
to each site (=lattice point). The edges between the lattice points
correspond to (nearest neighbour) interactions and  the energy of
a state $\sigma$ is the sum:
$$\sum_{\hbox{edges between lattice points}}E(\sigma_x,\sigma_y)$$
where in the sum $x$ and $y$ are the lattice points at each end of the edge,
and $E(i,j)$ (with $i$ and $j$ being $\uparrow $ or $\downarrow $) is the
local energy arising from the interaction along the edge.

The boundary conditions can be handled in many ways-one can wrap     
approximating rectangles on a torus creating periodic boundary 
conditions. Or one can simply neglect the interactions of the boundary
sites with neighbours outside $X$, (free boundary conditions), or one
may fix all the spins along the boundary according to some specified 
pattern (fixed boundary conditions). Since most of the contribution 
to the partition function will not involve the boundary, the 
asymptotic growth rate of the partition function should depend 
only on the whole system. This rate is called the free energy
per site:

$$ F= \lim_{N\rightarrow \infty} {1\over N}\log Z_X$$

This free energy may depend on several parameters. Certainly
the temperature is one of them, but there may be different 
horizontal and vertical interactions, an external field and
so on. 
\5
We will say a model is ``solved'' if $F$ is expressed as an
explicit function of its parameters. Given the complexity
of the function that may be involved in such a solution, one
may question the usefulness of a solution as opposed to 
the defining limit. But there are many cases in which the
explicit formula is simple enough to read off meaningful
results. There are also many other limits one might like
to calculate before saying the model is ``solved''.
\5
The most completely (non-trivial) solved model is the Ising
model in 2 dimensions. But we shall be more interested in 
another kind of model called a ``vertex model'' on a lattice,
where the state of the system is defined by assigning values
(in some indexing set) to the \underline{edges} of the lattice.
The ``ice-type'' model is a vertex model in which the indexing
set has two elements (corresponding to the presence or absence
of some kind of bond between neighbouring molecules) which can
be conveniently represented by arrows on the edges. Thus a 
state of an approximating square in a 2-d ice-type model 
might be as below:
\5
\begin{diagram}\label{ice}
\vpic{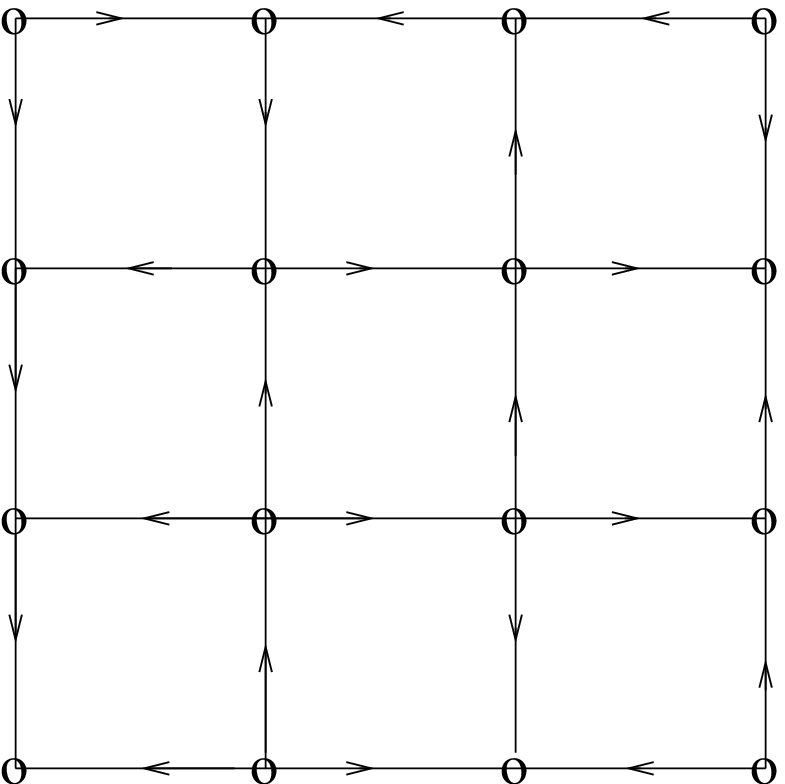} {2.4in}
\end{diagram}
\5
 The energy of a state of a vertex model is the sum of energy
contributions from each vertex. If the state is given each
vertex is surrounded by edges with indices on them so that
the energy is specified by assigning an energy to each
configuration of indices. In the ice-type model there
are 16 such configurations corresponding to 
the arrow configurations around a vertex.

The partition function is calculated using exponentiated
energies. The exponential of the energy is called the 
Boltzmann weight so that we have Boltzmann weights
$$R(a,b|c,d)$$ assigned to each local index configuration as below:
\5
\hspace{1in}\vpic{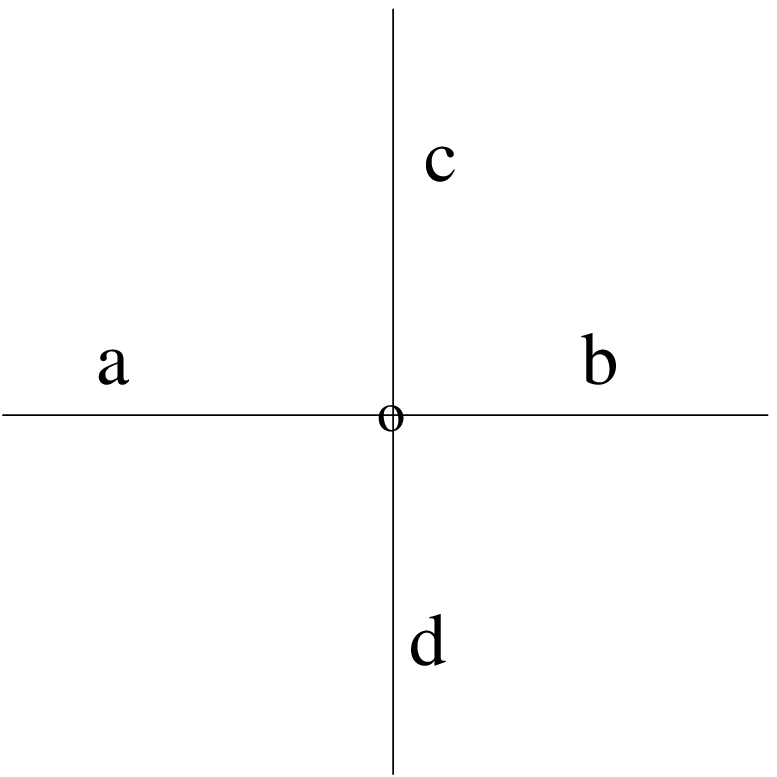} {1.5in}
\5
The partition function for the rectangular subregion $X$ is then
$$Z_X=\sum_{states}\prod_{vertices}R(a,b|c,d)$$

Where conventions must be adopted for how the indices surrounding
a vertex (in a given state) are to be used as indices in $R(a,b|c,d)$,
and the boundary conditions must be specified.
\begin{remark}
In a large part of the literature what we call ``$R$'' below is
called $\check R$ and $R=S\check R$ with $S$ as
in section \ref{decomposition}. We use our notation slightly reluctantly but
it seems that the more fundamental formalism is the one where 
our $R$-matrix is present but $S$ is not. And we do have the justification that
$R$ is the letter Baxter himself uses in \cite{baxter}.
In quantum group theory it is no doubt the other $R$ that is more
natural.
\end{remark}
\subsection{Transfer Matrices} 
Transfer matrices are a powerful method for translating
the problem of finding the partition function into a problem
of linear algebra. The basic idea is that the summation 
over indices in the partition function becomes the 
summation over indices in matrix multiplication. For instance
if one had a one dimensional vertex model with 
Boltzmann weights $R(a,b)$ the partition function for 
a lattice with $n$ sites as below (illustrated with $n=5$):
\5\5\5
\vpic{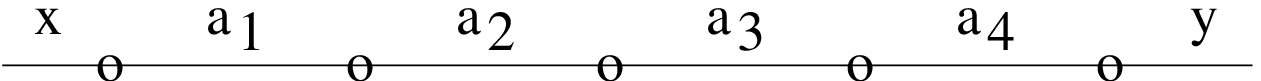} {3.5in}
\5\5\5
\noindent is readily seen to be the $(x,y)$ entry of the 
matrix $R^n$. The boundary conditions were fixed to be
$x$ at the left and $y$ at the right. If the boundary
conditions were periodic the partition function
would be $Trace(R^n)$.
\5
One is interested in the asymptotic behaviour as the 
subsystem $X$ tends to the whole infinite lattice and one
can use linear algebra techniques to understand the
asymptotic behaviour of $R^n$ (the behaviour is in general
governed by the largest eigenvalue.
We leave the solution of the  one dimensional
vertex model as an easy exercise.

\5
To apply the transfer matrix method to a two dimensional lattice
we simply think of each row of the lattice as being an atom
in a one dimensional lattice and construct its transfer matrix.
 The trouble is of course that
the the size of the transfer matrix will grow (exponentially) with
the size of the system.  And the boundary conditions will have to
be handled in a more complicated way. Let us first impose 
periodic horizontal boundary conditions. Then the transfer matrix
for a 2-d lattice built up from horizontal rows will be
$$T_{x_1,x_2,...,x_n}^{y_1,y_2,...y_n}=\sum_{a_1,a_2,...,a_n}
 R(a_n,a_1|x_1,y_1)R(a_1,a_2|x_2,y_2)....R(a_{n-1},a_n|x_n,y_n).$$ 
as explained diagrmattically below:
\5
\begin{diagram}\label{transferm}\vpic{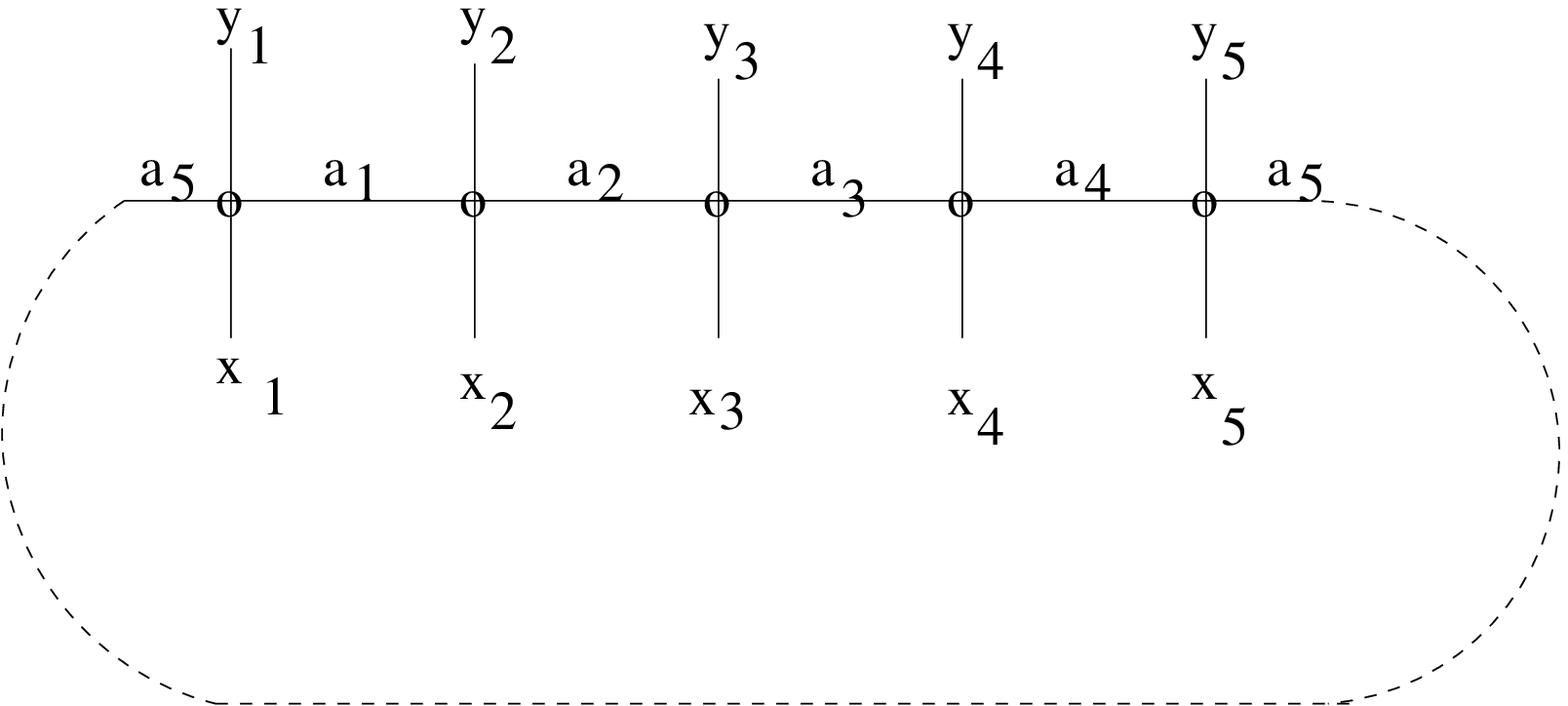} {3.5in}
\end{diagram}
\5
Because of the growing size of $T$, the problem of calculating
its largest eigenvalue becomes formidable and hopeless in general.
One of Baxter's great ideas was to look for models in which
the transfer matrices \underline{commute} \underline{with} \underline{each } \underline{other} for
different values of their parameters. Then they will have to 
have a common eigenvector and one may try to deduce enough
about how the eigenvalue depends on the parameter to determine
it completely. This part of the Baxter program - actual 
determination of the eigenvalues - has not been completely 
formalised, but a great machine has evolved for producing
examples of models with commuting transfer matrices. That
machine is \emph{QUANTUM GROUPS}.

\subsection{The Yang-Baxter equation.}
The diagram below illustrates what it means for 
the transfer matrix with value $\lambda$ (often called
the \emph{spectral parameter}) to commute with the
transfer matrix with value $\mu$ (periodic horizontal boundary conditions):

\vpic{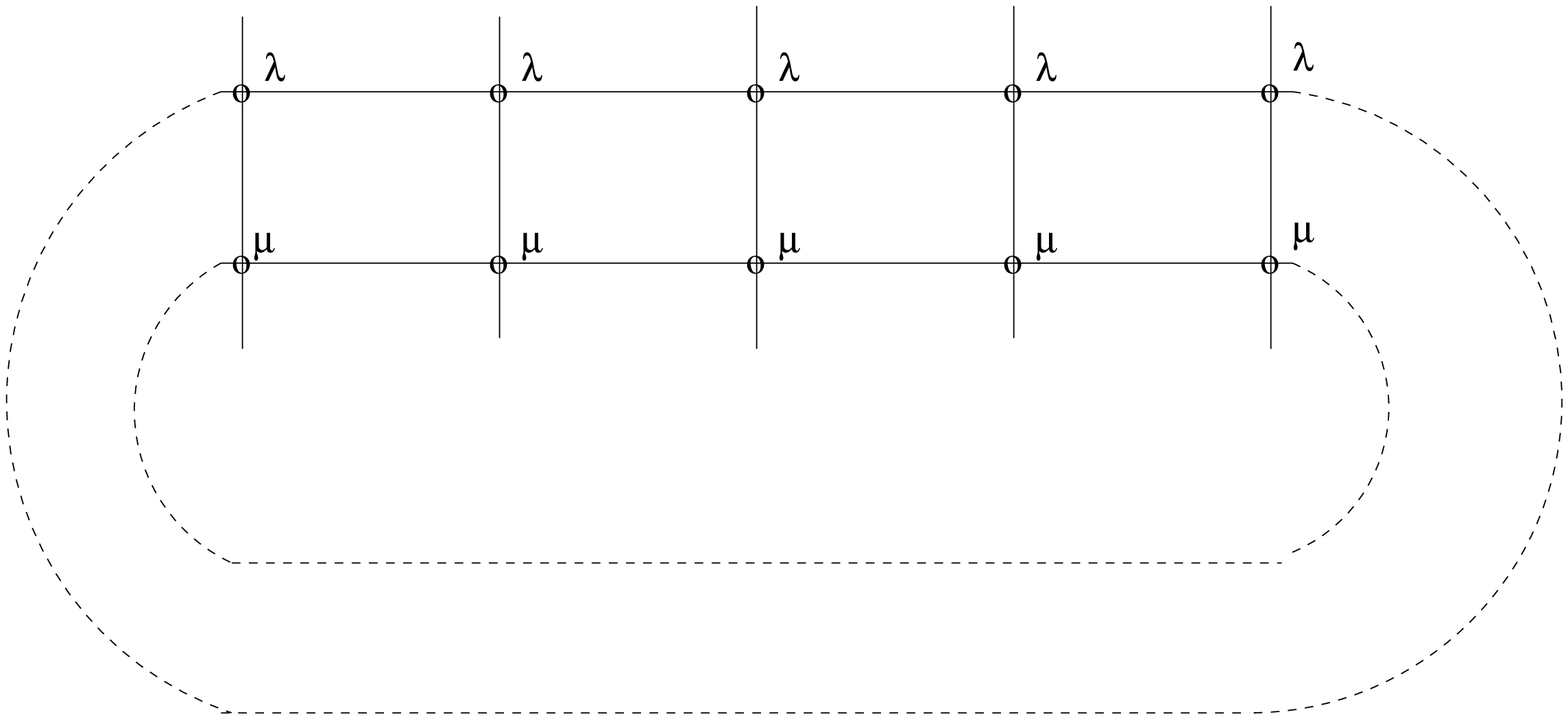} {4in}
\5
=\quad \vpic{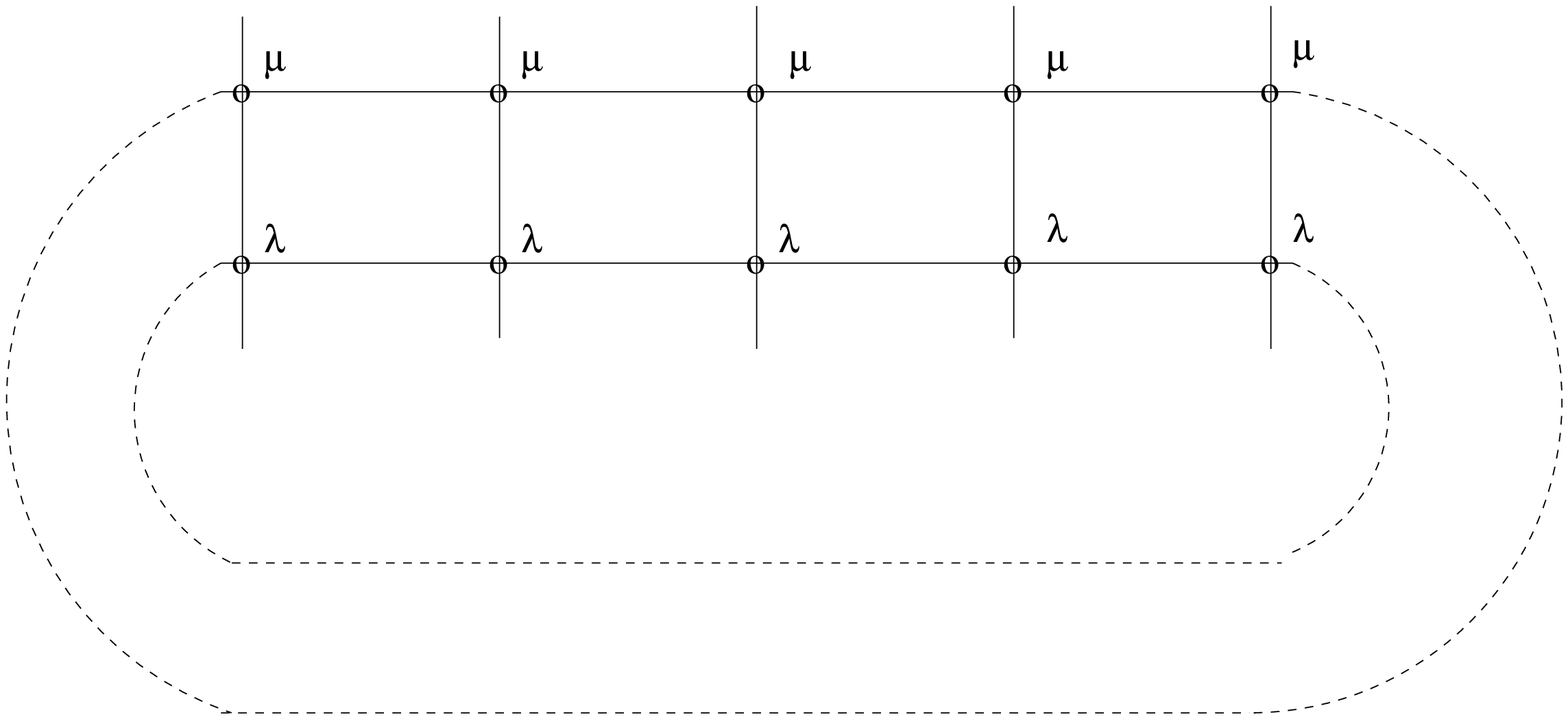} {4in}

Here we left out all 
indices, the convention being that indices are implicit on
the boundary edges and summed over for each internal edge.
And the value of the spectral parameter to be used for
the $R$ matrix is indicated near to the corresponding vertex
on the diagram.

If written out in full, the equations represented by the 
diagram form a huge system of highly non-linear equations for
the Boltzmann weights. The Yang-Baxter equation (YBE) is 
a set of equations involving ony 3 vertices which implies
that the transfer matrices commute.  With the same notational
conventions as above the YBE asserts the existence of a third value
$\rho$ of the spectral parameter (depending of course on
$\lambda$ and $\mu$) for which we have the following equation:
\5
\hspace{0.7 in}\vpic{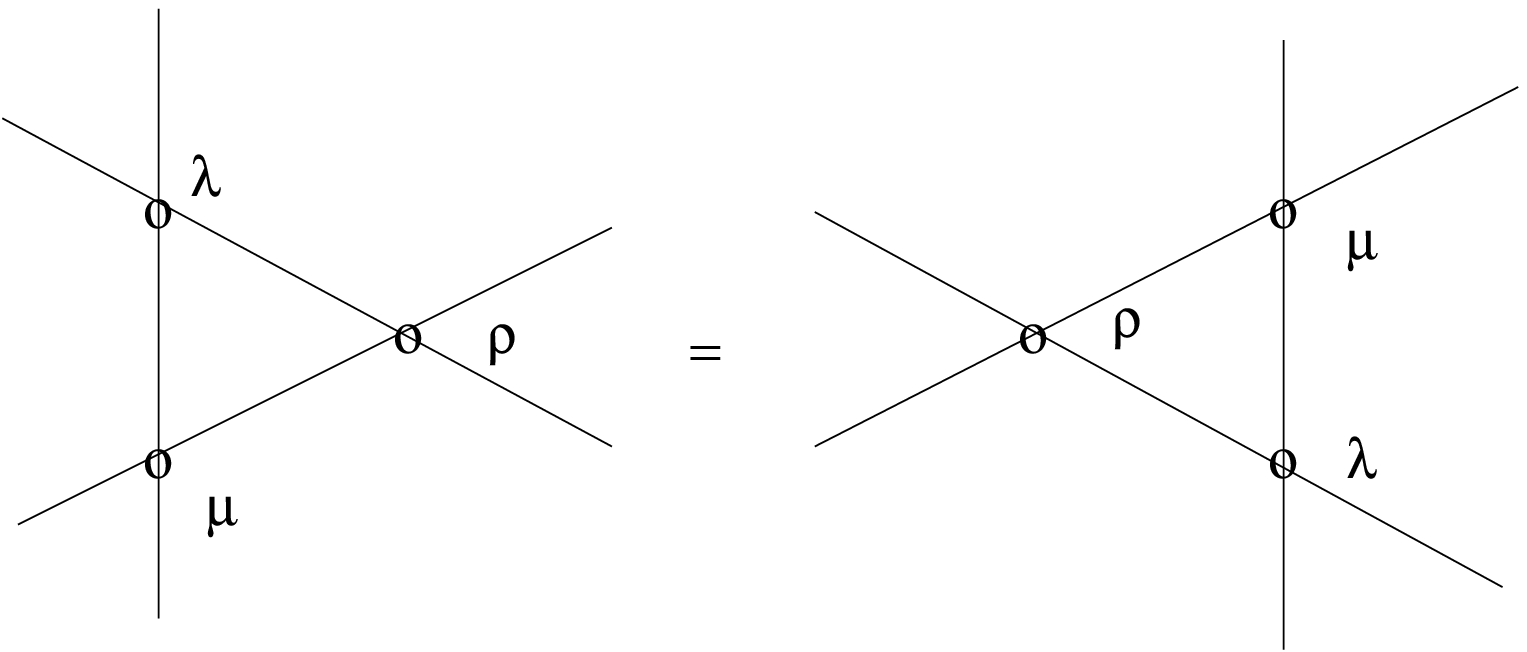} {3in}
\5

If we use $R(\lambda)$ to denote the matrix of Boltzmann weights
with parameter $\lambda$ then the YBE is, in the  notation of \ref{decomposition}:
\begin{formula} \label{ybe}
\quad$R_{12}(\lambda)R_{23}(\rho)R_{12}(\mu)=R_{23}(\mu)R_{12}(\rho)R_{23}(\lambda).$
\end{formula}

The argument that the YBE implies commuting transfer matrices is an 
elegant one which is entirely diagrammatic with our summation 
convention. We need to make the assumption that the  matrix of
Boltzmann weights for the third value $\rho$ is invertible. This
is precisely the condition that there is another R-matrix which
we will denote by the parameter ``$\rho^{-1}$'' for which:
\5\5
\vpic{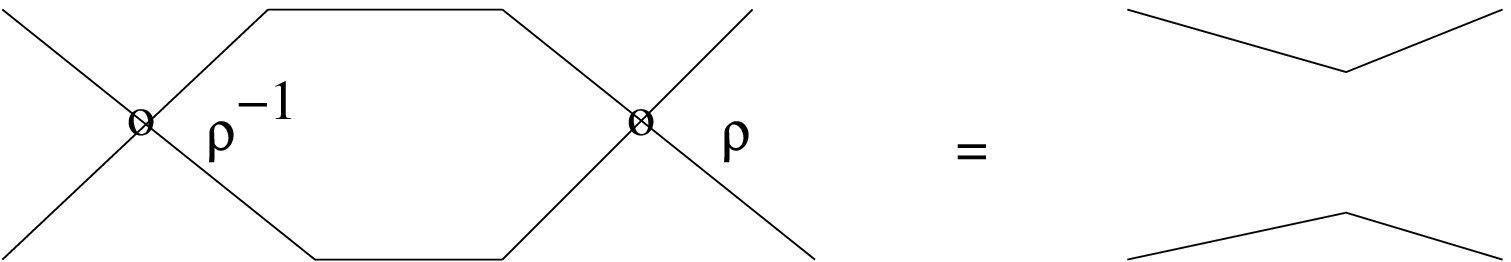} {3.5in}
\5
\noindent (which of course implies the same thing with $\rho$ and $\rho^{-1}$
interchanged). Note that it is rather important to associate the
correct indices of $R(a,b|c,d)$ to the correct edges of the the
diagram. How to do this will be obvious from the following argument
so we leave it to the reader.

Now take the picture representing one side of the equation for
commuting transfer matrices and insert the picture above
 for $\rho$ ``$\rho^{-1}$''
to obtain:
\5\5
\vpic{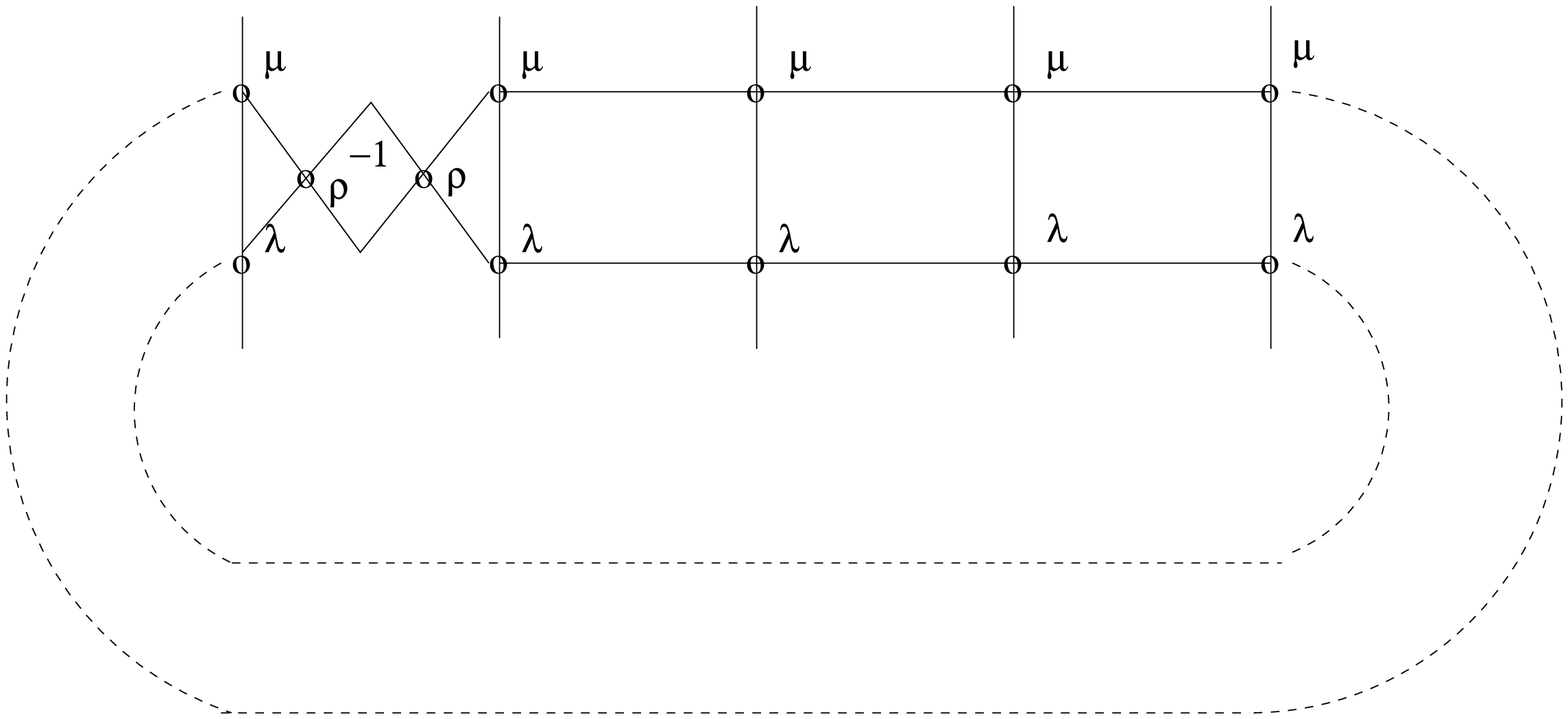} {4in}
\5

This does not change the partition function. Successive applications
of YBE move $\rho$ clockwise around the picture, swapping a $\lambda$
and a $\mu$ each time. After a few steps one obtains:
\5\5
\vpic{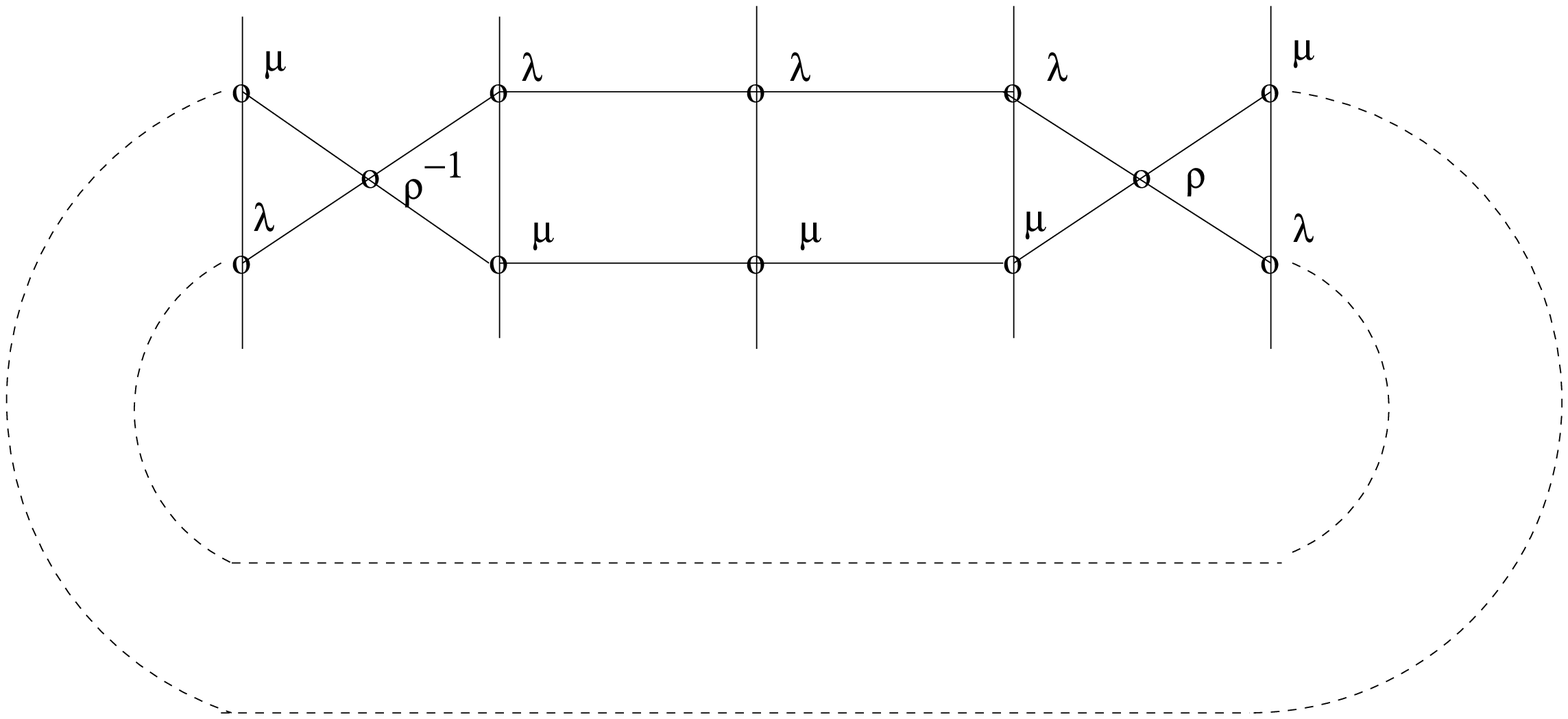} {4in}

\5
Eventually the $\rho$ comes right round the circle and meets its 
inverse with which it cancels and one obtains the other side of 
the commuting transfer matrix equation!!

\subsection{First solution of the YBE}

Apart from the most difficult question of how to proceed 
once we have commuting transfer matrices (which we do not
pursue here), the YBE raises many questions. One might
ask for instance if the $\rho$ is uniquely defined by
$\lambda$ and $\mu$ and if so what are the properties
of the operation $(\lambda,\mu)\mapsto \rho$.
But most of all there is the question of existence-
the YBE is a system of coupled cubic equations with
far more equations than unknowns. Without any discussion
at this point of how it was found, we present the following
$R-$matrix and claim it is a solution of YBE:

\begin{formula}\label{rmatrix}
 \hspace{1in}$R_q(x)=$

\5$${1\over xq-x^{-1}q^{-1}}
\left( \begin{array}{cccc}
xq^{-1}-x^{-1}q &0&0&0  \cr
0 & x^{-1}(q^{-1}-q)&x-x^{-1} &0 \cr
0&x-x^{-1}& x(q^{-1}-q)&0  \cr
0&0&0& xq^{-1}-x^{-1}q \cr
\end{array} \right)
$$
 \end{formula}
The assiduous reader may check directly that (suppressing $q$-dependence)
$$R_{12}(x)R_{23}(xy)R_{12}(y)=R_{23}(y)R_{12}(xy)R_{23}(x)$$
but we will see easier ways to check this later on.
Note that the factor $\displaystyle {1\over x^{-1}q^{-1}-xq}$ is arbitrary but
it yields the following pleasant properties:
\5
\5

\begin{description}
\item{(\romannumeral 1)}\quad$ R_1(x)=S$ \qquad \qquad(recall from \ref{decomposition}
that $S$ is the flip $v\otimes w\mapsto
w\otimes v$).
\item{(\romannumeral 2)}\quad $ R_q(1)=-1$.
\item{(\romannumeral 3)}\quad $R_q(x)^{-1}=R_q(x^{-1})$
\item{(\romannumeral 4)}\quad  If we define $R$ to be $\lim_{x\rightarrow 0} R_q(x)$
then

$\displaystyle {R=
\left( \begin{array}{cccc}
q^2 &0&0&0  \cr
0 &q^{2}-1&q &0 \cr
0&q& 0&0  \cr
0&0&0& q^{2} \cr
\end{array} \right)
}$

 satisfies the braid equation 
$$R_{12}R_{23}R_{12}=R_{23}R_{12}R_{23}.$$
\end{description}
\5 
Note that writing $x=e^\lambda$ and $q=e^\theta$ converts all
the entries of the matrix into hyperbolic sines. This R-matrix
is called a trigonometric solution of the YBE.
\5
Since the entries of the matrix are supposed to be Boltzmann weights,
to be of interest to statistical mechanics there must be values of
$x$ and $q$ for which the entries of the matirx are all non-negative.
The global multiplying factor is neither here nor there so 
we see that it suffices to choose $x$ and $q$ positive with $x>x^{-1}$
and $q^{-1}>q$.

\subsection{The Ice-type model and the Potts model.}
 Recall that a Boltzmann weight of $0$ corresponds to 
infinite energy-i.e. a forbidden configuration. If we
look at the positions of the zeros of our R-matrix \ref{rmatrix}
and think of the rows and columns as indexing arrow configurations
as in \ref{ice} we see that the configurations allowed by the 
R-matrix are the following:
\5
\vpic{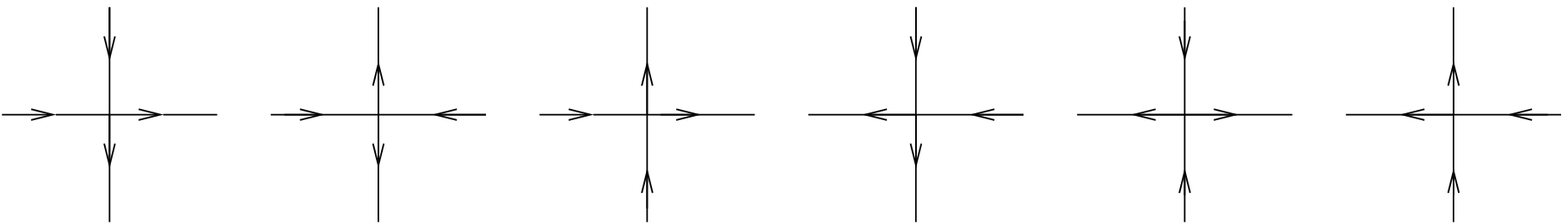} {4.5in}
\5
Excluding all but these particular configurations has something to do
with the particular physical system under configuration.
Lieb solved the ice-type model, not by the method
of commuting transfer matrices but by the so-called Bethe Ansatz.
(For details on all this see \cite{baxter}.)

\5
The Potts model is not a vertex model. It is like the Ising model
in that individual ``spins'' are located on the vertices of the
lattice and a state of the system is specified by assigning a
``spin'' value from $1$ to $Q$ to each of the vertices. The interactions occur
along the edges of the lattice so that the total energy of 
a state is 
$$\sum_{\hbox{edges of the lattice}} E(\sigma, \sigma')$$
where we have suppressed the approximating rectangle and $\sigma$ and
$\sigma'$ denote the spin values at the ends of the particular edge
being summed over. Thus the partition function is
$$\sum_{states}\hspace{7pt} \prod_{edges}w(\sigma,\sigma')$$
where the Boltzmann weights are the exponentiated energies as usual.
Thus from a purely mathematical point of view the only data 
for the lattice model is the $Q\times Q$ matrix $w(\sigma,\sigma')$
of Boltzmann weights. If the edges of the lattice are not directed
this must be a symmetric matrix though the geometry of the lattice
may allow, say, different Boltzmann weights for vertical or horizontal
interactions.

\5
The Potts model is defined by the property that the spin states have
no structure other than being different so that the Boltzmann weight
$w(\sigma,\sigma')$ depends only on whether $\sigma=\sigma'$ or not.
If $V=\mathbb C^Q$ with usual basis $v_\sigma$,
then the transfer matrix which creates a new row
with $n$ spins of the lattice will be a linear map from $\otimes^n V$
to itself. To organise the transfer matrix we introduce the maps
$p:V\rightarrow V$ with all matrix entries equal to $1\over \sqrt Q$, and
the map $d:V\otimes V\rightarrow V\otimes V$ with 
$d(v_\sigma \otimes v_{\sigma'})= \delta_{\sigma, \sigma'}v_\sigma\otimes v_\sigma$.
We then put $E_{2i-1}=1\otimes 1\otimes...\otimes p \otimes 1... \otimes 1$ with
the $p$ in the $i$th. tensor position, and 
$E_{2i}=\sqrt{Q} d_{i(i+1)}$ using the notation of section \ref{decomposition}. 
Then it is an easy exercise to show that, for the Potts model, 
the transfer matrix with free horizontal boundary conditions is
a multiple of
\begin{formula} \label{transfermatrix}
\qquad \qquad $\displaystyle{\prod _{i=1}^{n-1} (aE_{2i}+1) \prod _{i=1}^n(bE_{2i-1}+1)}$
\end{formula} 
where $a$ and $b$ are determined by the horizontal and vertical Boltzmann
weights respectively (note for instance that necessarily, up to a constant,
$w=Ap +1$ where $w$ is the matrix given by the Boltzmann weights).
\5
You may be wondering about the bizarre normalisations we have used in
defining $p$ and $d$ and the strange indexing of the $E_i$'s.
The reason is to get the nicest possible algebra going. It should
be checked that the $E_i$ satisfy the following relations:

\begin{description}
\item{\begin{formula}\label{tl1}\qquad $\displaystyle{ E_i^2=\sqrt{Q}E_i}$ \end{formula}}
\item{\begin{formula}\label{tl2}\qquad $\displaystyle{ E_i E_{i\pm 1}E_i =E_i}$ \end{formula}}
\item{\begin{formula}\label{tl3}\qquad $\displaystyle{ E_iE_j=E_jE_i}\qquad \hbox{if}\quad |i-j|\geq 2$ \end{formula}}
\end{description}

These relations are known as the Temperley-Lieb relations and are somewhat
magical. It is fun to check that $P=E_1E_3E_5...E_{2n-1}$
has the property that 
$$PxP=\phi(x)P$$ where $x$ is any word on the $E_i$'s and $\phi(x)$ is
a real valued function of $x$. This means that if $X$ is any element
in the algebra generated by the $E_i$'s then $PXP=\phi(x)P$ for
some linear functional $\phi$ on that algebra. Moreover in terms of
our statistical mechanical model this functional $\phi$ gives 
precisely the partition function for free vertical boundary conditions!
Thus in principle the partition function for a rectangular lattice with
free boundary conditions is entirely determined by the Temperley-Lieb relations.

\5
So what? To answer this we return to our $R$-matrix \ref{rmatrix}
for the ice-type model.
Put 
\begin{formula}\label{modelE}

\5
$\displaystyle{E=\qquad
\left( \begin{array}{cccc}
0 &0&0&0  \cr
0 & q^{-1}&1 &0 \cr
0&1& q&0  \cr
0&0&0&0 \cr
\end{array} \right)
}$

 \end{formula}

Then two things are true. First, if we define $E_i$ on $\otimes^n\mathbb C^2$
as $E_{i(i+1)}$ with the notation of \ref{decomposition}, then the Temperley-
Lieb relations hold, and second, $R_q(x)$ is a linear combination of
$E$ and the identity. It is thus not surprising that, with the appropriate
boundary conditions, the partition function for the ice-type model is the
same as that of the Potts  model with a (physically bizarre) change of variables. In fact
this is only true if the horizontal and vertical interactions of the
Potts model satisfy the relation $a=b$, known as ``criticality'' for
various reasons. In \cite{templieb}, Temperley and Lieb showed the equivalence of
the ice-type model and the critical Potts model on a square lattice using
the relations \ref{tl1}\ref{tl2} and \ref{tl3}. This equivalence 
was subsequently understood on a general planar graph (for the Potts model)
and its "medial graph" (for the ice-type model). For a beautiful treatment
see chapter 12 of \cite{baxter}.

\subsection{How to remember the formula.}

My personal way of reconstructing the formula \ref{rmatrix} from 
simpler ones involves the Hecke algebra of type $A_n$. This is
the algebra with generators $g_1, g_2,...,g_{n-1}$ and relations

\begin{description}
\item{(h1)}\qquad $g_i^2=(q-1)g_i +q id$
\item{(h2)}\qquad $g_i g_{i+1} g_i = g_{i+1}g_ig_{i+1}$
\item{(h3)}\qquad $g_ig_j=g_jg_i$\quad if  \quad $|i-j|\geq 2$.
\end{description}

Here I am faithfully reproducing a constant disagreement in the
literature over the meaning of $q$. In our Hecke algebra relations
we are using $q$ as in \cite{bourbakilie}, which is natural in its context as 
the number of elements in a finite field. The $q$ in \ref{rmatrix}
is the square root of this $q$.

The relations $h2$ and $h3$ are the braid relations which we
have seen as the limit of the YBE as the spectral parameter 
tends to infinity. In this Hecke algebra case we can reconstruct the YBE from
the braid relations as follows:
\5
Step 1: Renormalize the $g_i$'s as $G_i$'s so that relation $h1$ 
becomes 
$$G_i+G_i^{-1}= k\hspace{3pt} id$$

Step 2: Define $$ R_i(x) = xG_i+x^{-1}G_i^{-1}$$.

Then it is an exercise to prove that, in the presence of the braid
relations the YBE is equivalent to

\begin{formula}\label{baxtercond}
$\displaystyle \qquad G_1G_2^{-1}G_1 + G_1^{-1}G_2G_1^{-1} =G_2G_1^{-1}G_2 + G_2^{-1}G_1G_2^{-1}$
\end{formula}

It is immediate to show \ref{baxtercond} from $G_i+G_i^{-1}= k\hspace{3pt} id$.
I do not know of any other solutions to \ref{baxtercond}.

Of course this begs the question of how to get appropriate solutions
of the Hecke algebra relations. One way is to obtain $g_i$'s from
the Temperley-Lieb $E_i$'s by 
$g_i=qE_i-1$ (which $q$ is this??) and this is indeed how I put together
\ref{rmatrix}. But there are other solutions as we shall see.

\section{Local Hamiltonians}
 A \emph{quantum spin chain} is a one dimensional array of 
spins. If the Hilbert space describing an individual spin
is $\mathbb C^2$, then a quantum spin chain with $n$ spins
will be described by $\otimes^n \mathbb C^2$. If $T$ is a Hermitian
$2\times 2$ matrix defining an
observable for a single spin, then that observable for the $k$th. spin
in the chain is (with $T$ in the kth. slot):
$$T_k=1\otimes 1\otimes... \otimes T\otimes 1\otimes...\otimes 1.$$
If $H$ is the Hamiltonian for the spin chain, the observables $A$
evolve according to $A_t = e^{iHt}Ae^{-iHt}$.
By writing down the correlations between observables of spin $k$ at (discretized)
time $t$ one sees a strong similarity with expected values of 
spins in the Ising model whose $x$ coordinate is given by $k$ and
$y$ coordinate by the time $t$, provided one takes as Hamiltonian
the logarithm of the transfer matrix (times $\sqrt{-1}$). This is
generalised into a powerful equivalence between 2-dimensional statistical
mechanics and 1-dimensional quantum mechanics provided time
can be anaylytically continued to imaginary time.

This suggests that perhaps the transfer matrices of statistical mechanical
models can be used to create Hamiltonians for quantum spin chains. In order
to satisfy locality conditions, a Hamiltonian should be expressible as a
sum of terms each one only involving spins close to each other on the 
lattice. The simplest would be a nearest neighbour interaction and 
if it is translation invariant it must be of the form
$\sum_i H_{i(i+1)}$ in the notation of \ref{decomposition}
where $H$ is some self-adjoint operator on $\mathbb C^2\otimes \mathbb C^2$.
An ingenious way to do this is to take the logarithmic derivative of the transfer matrix with respect
to the spectral parameter at some appropriate value of the spectral
parameter. By the conditions after \ref{rmatrix} it is clear that the
right value is $x=1$.

Since the transfer matrix is multilinear in the $R$ matrices used
at the vertices we see that the derivate with respect to $\lambda$
of
\5
\hspace {1.0in}\vpic{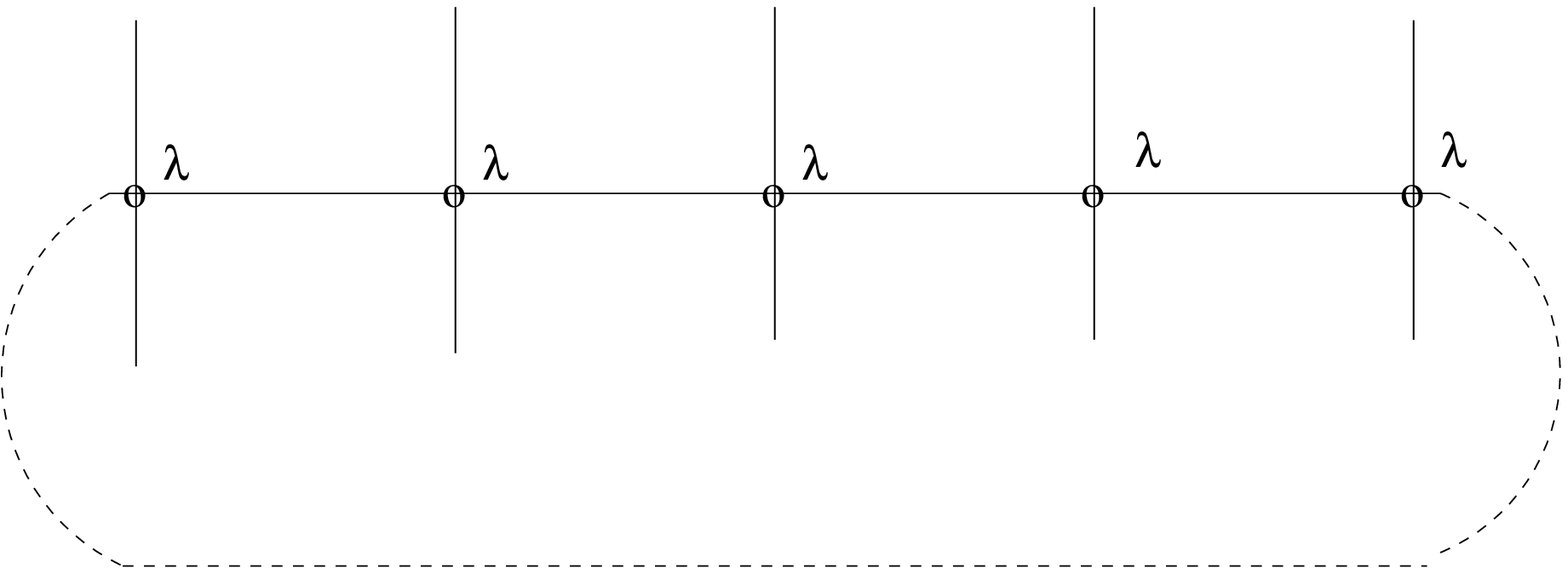} {3in}
\5
\noindent is the sum over all ways of putting in one
$\lambda'$ of 
\5
\hspace {1.0in}\vpic{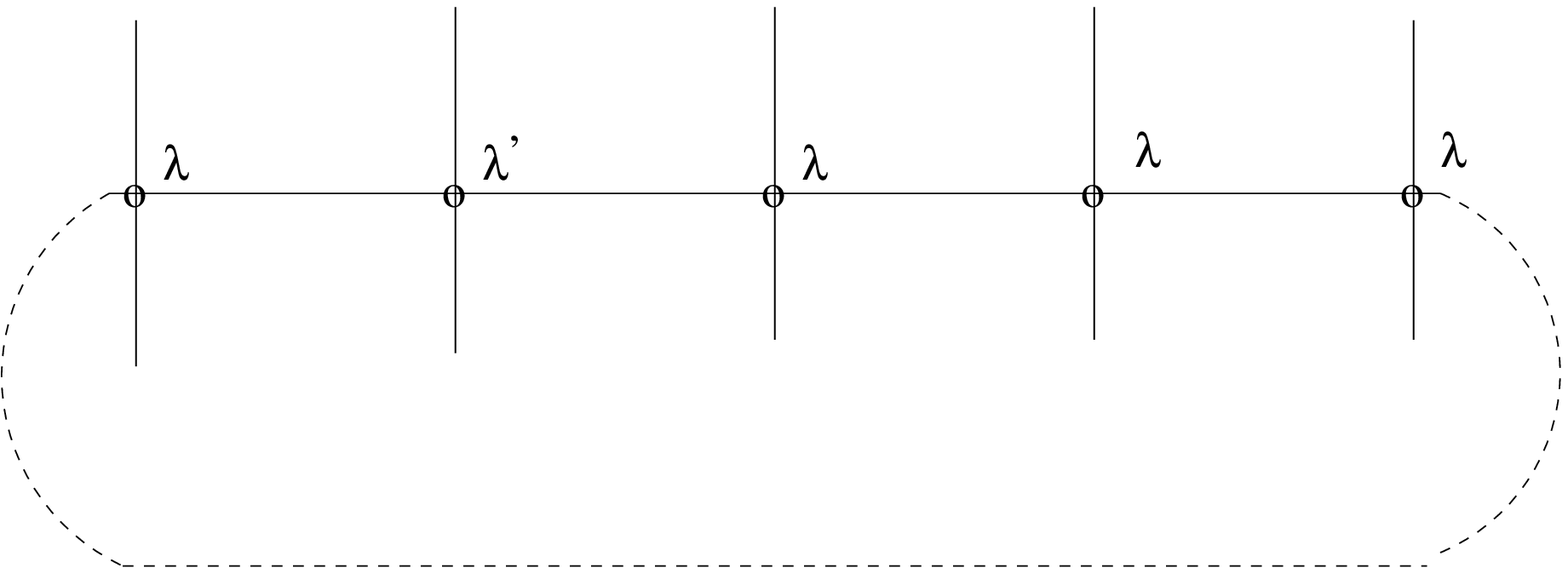} {3in}
\5
\noindent where we have symbolically used $\lambda'$ to stand for the derivative
of the $R$-matrix with respect to $\lambda$. If we use \ref{rmatrix},
the sign of $R_q(1)$ is irrelevant so we  see that
the transfer matrix \ref{transferm} is represented diagrammatically by:
\5\5
\hspace {1.0in}\vpic{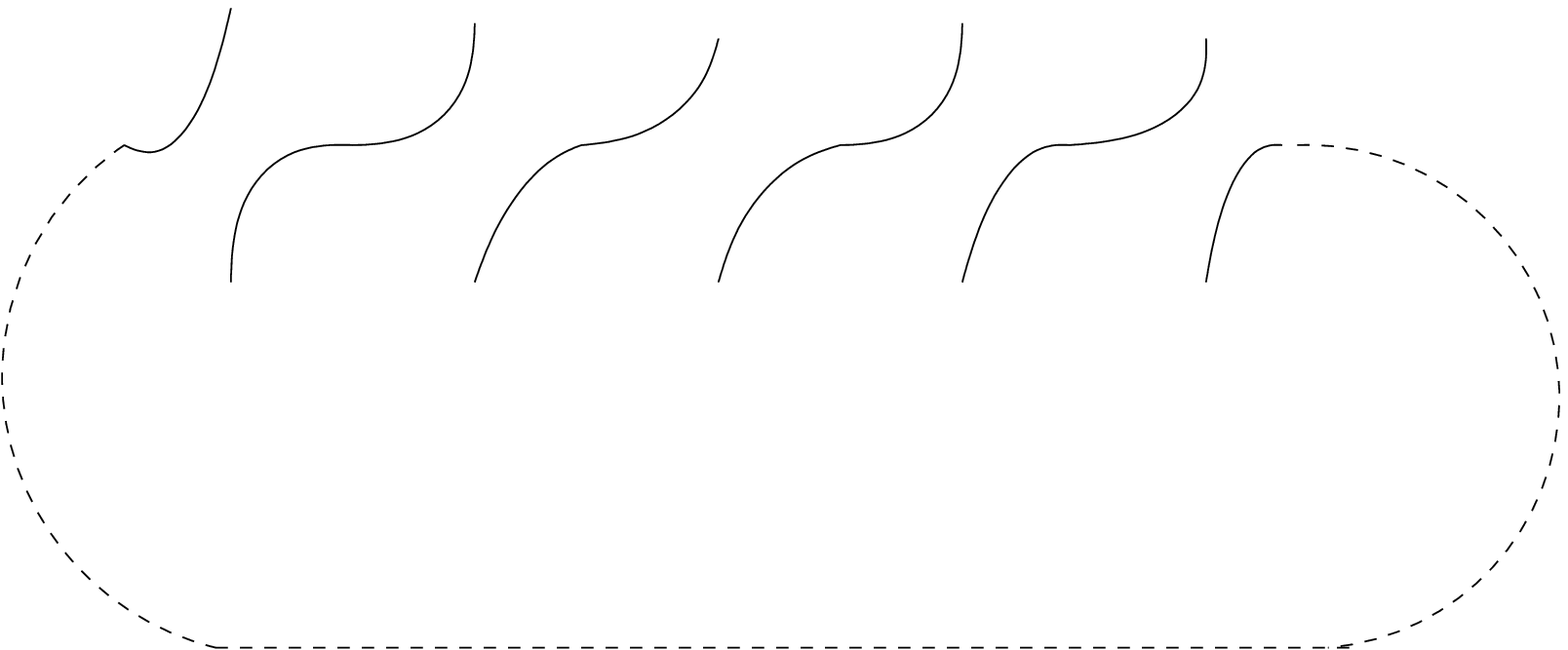} {3in}
\5
which is of course a rotation if it is represented in a cylinder.
\5

So if we write $T(x)$ for the transfer matrix using \ref{rmatrix},
a typical term in the logartithmic derivative $$\displaystyle{T(1)^{-1}{dT\over dx}|_{x=1}}$$ 
can be seen by multiplying the diagrams  as below 
with $R'(0)$ being the derivative of \ref{rmatrix}
with respect to $x$ at $x=1$: 
\5
\hspace{1in} \vpic{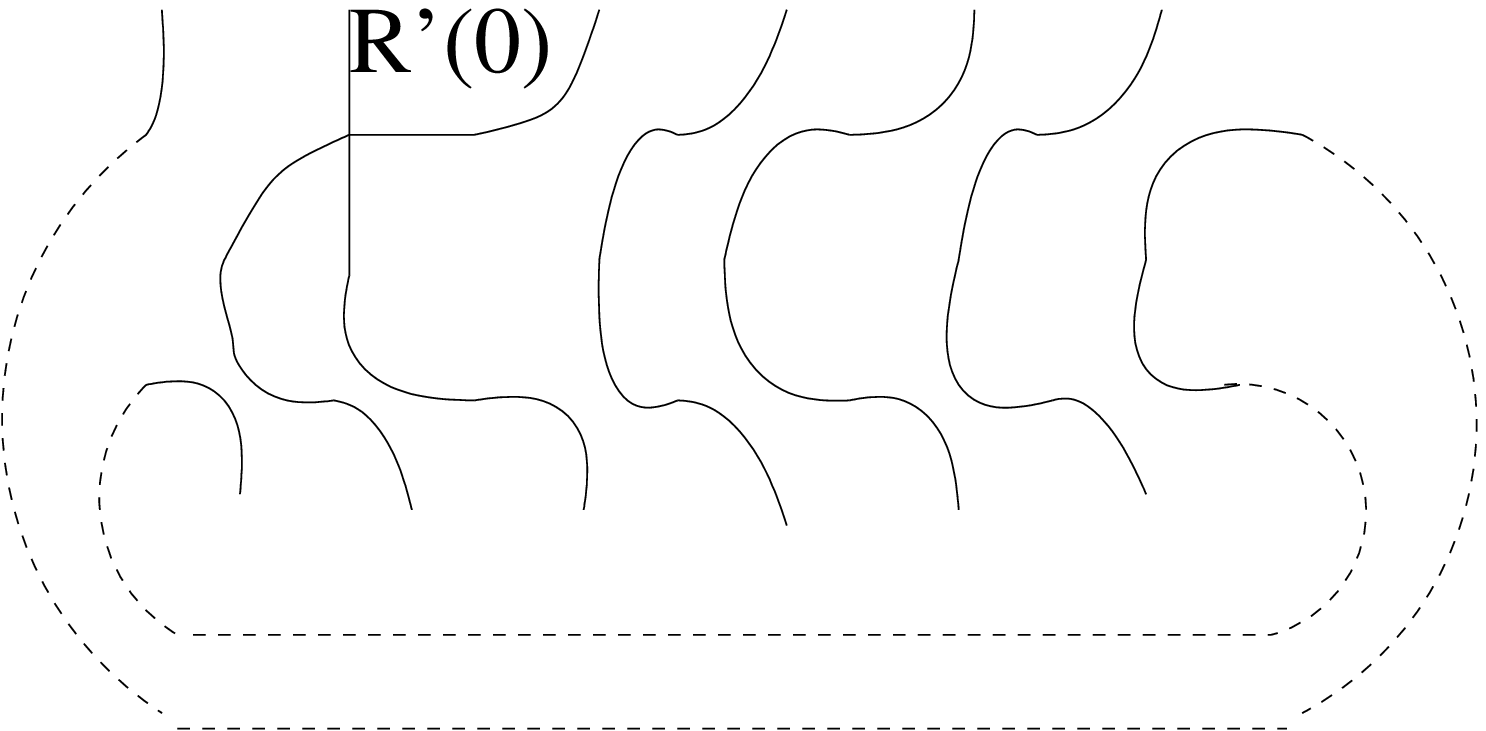} {3in}
\5

On cleaning up the picture we see that the logarithmic derivative is
the sum of all matrices whose diagrams are as below, the crossings
occuring between the $i$th. and $(i+1)$th. strings from the left,
with periodic boundary conditions so that the last term would involve
a crossing between the first and last strings: 
\5
\hspace{1in}\vpic{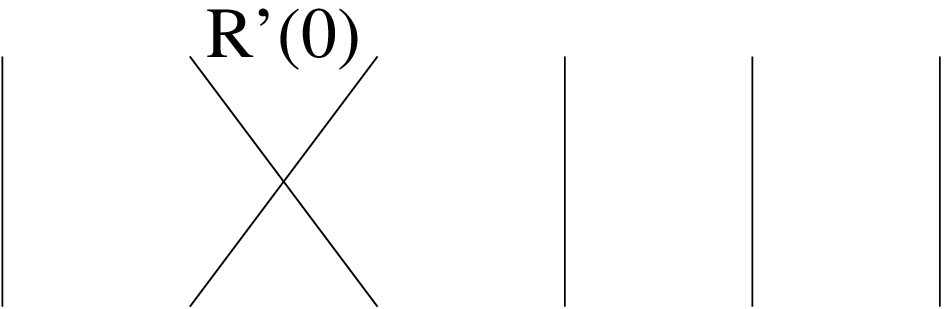} {3in}
\5
This sum of matrices is clearly a Hamiltonian with nearest neighbour
interactions provided it is positive self-adjoint. This positive 
condition on the $R$-matrix may be quite different from the positivity
of the Boltzmann weights so $R$-matrices which do not work in the
statistical mechanics world may work here. The moral is that if
you find a solution to the Yang-Baxter equation that does not
admit positive Boltzmann weights, don't necessarily condemn it to
the trash as it may give you a solvable quantum spin chain. 
\5 
  I have made the arguments for commuting transfer matrix and
local Hamiltonian completely diagrammatic so they will work in
any situation where the diagrams make sense as multilinear 
maps on their inputs. For instance the arguments work just as
well for the Potts model and other models known as ``IRF models''.
Several people have axiomatised the diagram calculus-see \cite{kuperberg},\cite{bawest}.
I have developed a specific formalism, called ``planar algebras'' whose special 
features were driven by subfactors.
Certainly the arguments for commuting transfer matrices and
local Hamiltonians work in planar algebras.
\5 
We should not forget what we have achieved with  the local Hamiltonian-
since the transfer matrices all commute among themselves, they 
commute also with the local Hamiltonian so we are armed with
a large family of operators commuting with the time evolution
which should be extremely useful in diagonalising it. If, for
physical reasons, one wanted \emph{local} expressions for these
constants of the motion one could take higher logarithmic derivatives
of the transfer matrices with respect to the spectral parameter.
\5
In the particular case of our $R$-matrix \ref{rmatrix} of we do
the computation we find the local Hamiltonian
$$\sum_{i=1}^n H_{i(i+1)}$$ where the indices are taken mod $n$ and,
 up to addition of some constant
matrix which only changes the whole Hamiltonian by addition of a constant,

$$H_{i(i+1)}=\left( \begin{array}{cccc}
q^{-1}+q &0&0&0  \cr
0 & -(q^{-1}-q)&2 &0 \cr
0&2& q^{-1}-q&0  \cr
0&0&0& q^{-1}+q \cr
\end{array} \right)
$$
which may mean more to physicists written in terms of the Pauli spin matrices:

$$H_{i(i+1)}=\sigma_x \otimes \sigma_x +\sigma_y\otimes \sigma_y +{1\over 2}
\Bigg((q+q^{-1})(id+\sigma_z\otimes \sigma_z) +
 (q-q^{-1})(\sigma_z\otimes id -id\otimes \sigma_z)\Bigg)
$$

The presence of the term multiplied by $q-q^{-1}$ is the only difference 
between this and what is known as the XXZ Hamiltonian. But in fact as
Barry McCoy pointed out to me these terms cancel when one performs 
the sum over $i$. So the XXZ Hamiltonian admits a large family of explicit
matrices that commute with it.

\subsection{Spin Models}
As we have mentioned, quantum groups gave a machine
for creating large families of solutions of the Yang-Baxter equation and hence 
statistical mechanical models with commuting transfer matrices, and
quantum spin chains with many commuting Hamiltonians. But the $R$-matrix
coming from a quantum group is that of a vertex model. I have long
wondered if there is any such machine that would produce what
I call ``spin models'', that is to say generaisations of the Potts
model with an arbitrary (symmetric) matrix of Boltzmann weights.
There does not seem to be such a machine.  Searches for such models
have been very combinatorial and though they have led to some insights
in combinatorics (see \cite{bbj}) there have been few new statistical mechanical
models. The one spectacular new spin model was discovered by
Jaeger in \cite{jaeger}. He only gives the knot theoretic solution of the
Yang-Baxter equation but it can be easily ``Baxterised'' as in
\cite{jones:fz} to give a solvable model.

The idea that led to the Jaeger model was to look for models
whose matrix $w(\sigma,\sigma')$ 
of Boltzmann weights was the next simplest after
the Potts model. In the Potts model this matrix only has two
different entries, one on the diagonal and one off. The first
generalisation of this would be to matrices with three distinct
entries-one on the diagonal and two others, say $x$ and $y$.
One may then construct a graph whose vertices are the indices
of the matrix entries (i.e. the spin states per site) with
an edge connecting $\sigma$ and $\sigma'$ if
$w(\sigma,\sigma')=x$. This puts one in the world of association
schemes and their Bose graphs. By applying this kind of idea
Jaeger found two remarkable things:

(a) That there is a solvable model as above for which the
underlying graph is the Higman-Sims graph on $100$ vertices!
(The automorphism group of this graph is the Higman-Sims group,
one of the first sporadic finite simple groups.)

(b)Together with a couple of simple examples, the Higman Sims graph is 
the \underline{only} known graph that can work.

At this stage no-one has gone on to solve the Jaeger model.
There is a Temperley-Lieb like duality with a quantum group
R-matrix so that the bulk free energy is not too interesting, but
the correlation functions must support representations of the
Higman-Sims group.

\section{Subfactors.}

\subsection{Factors}
For a complete change of pace we treat a topic in analysis. Von Neumann
algebras are self-adjoint algebras of bounded linear operators on 
Hilbert space which contain the identity operator and
are closed under the topology of pointwise convergence
on the Hilbert space. Factors are von Neumann algebras with trivial
centre. We do not want to go into the details more than that but
we can suggest the notes from a course given by the author, accessible from his
home page, for anyone who wants to know more.
Thus a subfactor is a pair $N\subseteq M$ of factors.
To avoid technical difficulties we will only talk about the case of
``type II$_1$'' factors which are those which are infinite dimensional
but have a trace $tr$ which is a linear functional satisfying
$$tr(ab)=tr(ba)$$ and can be normalised so that $tr(1)=1$, in which
case $$tr(x^*x)>0 \hbox{    for   } x\neq 0.$$

\subsection{Index}
There is a notion of index for a subfactor, written $[M:N]$. If
$[M:N]<4$ it was shown in \cite{jones:index} that it must be one of the numbers
$4\cos^2\pi/n$ for $n=3,4,5,...$. The key ingredient in the proof
of this result was the construction of a tower of factors from 
the original pair and certain operators satisfying the Temperley-Lieb relations!
To be precise one can construct an orthogonal projection $e_N$
from $M$ to $N$ defined by the formula 
$$tr(xe_N(y))=tr(xy)$$ and then show that if $[M:N]<\infty$,
the algebra $<M,e_N>$ of linear operators on $M$ generated
by $M$ (acting by left multiplication) and $e_N$ is again a II$_1$ 
factor and $[<M,e_N>:M]=[M:N]$.
 This constructs the beginning of the tower:$$N\subseteq M\subseteq M_1= <M,e_N>.$$
To continue just repeat the construction to obtain $M_{i+1}=<M_i,e_{M_{i-1}}>$.

If we write $e_i=e_{M_i}$ for short and renormalise by $E_i=\sqrt{[M:N]}\hspace{2pt}e_i$
then relations \ref{tl1},\ref{tl2} and \ref{tl3} all hold. Moreover
there is *-structure for which $E_i^*=E_i$ and the trace has to 
satisfy $$tr(x^*x)>0 \hbox{    for   } x\neq 0.$$ A careful examination
of this property of the trace on the algebra generated by the $E_i$'s
proves the result about the ``quantized'' index values-see \cite{ghj}.

To see a more compelling similarity with the previous sections, note
that the tower $M_i$ can also be constructed as 
$$M_i= M\otimes_N M\otimes_N M \otimes_N...\otimes_NM$$ 
with $i+1$ copies of $M$ in the tensor product. 
(For anyone who has not seen the tensor product over non-commutative 
algebras, $M\otimes_N M$ is the quotient of the vector space tensor
product $M\otimes M$ by the subspace spanned by 
$\{xn\otimes y-x\otimes ny|x,y \in M \hbox{ and } n\in N\}$.)
The $E_i$'s can then be constructed with the same asymmetry between
odd and even $i$ as in the Potts model. Thus 
$$E_1(x\otimes y\otimes z...)=\sqrt{[M:N]}e_N(x)\otimes y\otimes z...$$
$$E_2(x\otimes y\otimes z...)=\sum_i xy\lambda_i\otimes \lambda_i^*\otimes z...$$
and so on, where $\lambda_i$ is any ``orthonormal basis'' for $M$ over $N$,
i.e. $\sum_i\lambda_i e_N \lambda_i^* = 1$ in $<M,e_N>$. The existence of
the $\lambda_i$ is an easy consequence of the original work of Murray and
von Neumann-see \cite{pp}. 
It is easy to show that the $E_i$'s acting on
$M_i= M\otimes_N M\otimes_N M  \otimes_N...\otimes_NM$ satisfy the TL relations
\ref{tl1},\ref{tl2} and \ref{tl3}.

In their first paper on subfactors, Pimsner and Popa discovered precisely
the representation of the TL relations given by \ref{modelE}! Not long 
afterwards D. Evans noticed the connection with statistical mechanical
models.

\subsection{A speculation on fusion/entanglement/interaction.}

I never fail to be struck by the relationship between the tensor 
product $M\otimes_N M$ and two interacting spins on the spin chain.
Connes had earlier introduced the notion of a ``correspondence''
between von Neumann algebras $M_1$ and $M_2$ which is a Hilbert 
space $M_1-M_2$ bimodule and had defined a surprisingly subtle
notion of tensor product of bimodules-see \cite{connes:ncg} and \cite{sauvageot}. Perhaps the 
kinematics for a system of two interacting quantum systems with
Hilbert spaces ${\cal H}_1$ and ${\cal H}_2$ with some common 
observables given by a von Neumann algebra $M$ is the Connes tensor
product ${\cal H}_1 \otimes_M {\cal H}_2$?

If two quantum systems described by Hilbert spaces
${\cal H}_1$ and ${\cal H}_2$ were so entangled that an 
an observable $x$ on one were equivalent to an observable
$y$ on the other then if $\xi \in {\cal H}_1$ and $\eta \in {\cal H}_2$
are vectors defining states, there should be no difference between
$x\xi\otimes \eta$ and $\xi \otimes y\eta$. If moreover the identification
of the $x$'s with the $y$'s were implemented by an antiisomorphism
$\phi :M\rightarrow \phi(M)$ from some von Neumann algbebra $M$ on
${\cal H}_2$ to a von Neumann algebra on ${\cal H}_1$ we would have a 
right action of $M$ on ${\cal H}_1$ and there would be no difference
between $\xi x\otimes \eta$ and $\xi \otimes x\eta$. We would be
forced to take the Connes tensor product ${\cal H}_1 \otimes {\cal H}_2$.

The behaviour of the Connes tensor product is very rich and contains
the theory of subfactors. Thus one could account for physicists'
assertions that the Hilbert space which describes several Chern
Simons ``particles'' is not the tensor product but a more 
complicated structure.  Thus this notion of fusion might be relevant
for the systems proposed by Freedman et al. in \cite{freedman}
in connection with quantum computing. 
The approach of Wassermann in \cite{wassermann} to the
fusion of loop group representations fits exactly into our 
framework and produces the right fusion algebra. 

Note that this notion of fused systems is much stronger than
the usual notion of entanglement where individual states
of a combined system may be more or less entangled.

The Hilbert spaces of the Andrews-Baxter-Forester and other IRF 
models for $n$ sites on a lattice furnish another example.
The natural basis for these Hilbert spaces is a basis of paths
and their dimensions are not simply powers of a given integer.
This could be explained if one supposes that proximity on the
lattice causes a large algebra of observables for one particle
to be identified with observables for its neighbour.

The Connes tensor product is easy to describe in finite dimensinos.
If the algebra $M$ is the $n\times n$ matrices and it acts (unitally) on the
right on a finite dimensional vector space $V$, $V$ may be 
identified tiwh the $p\times n$ matrices for some $p$, the right
action being matrix multiplication. Similarly a left $M$-module
structure on $W$ means that $W$ is isomorphic to the $n\times q$
matrices for some $q$. The tensor product $V\otimes_M W$ is then
the $p\times q$ matrices. Direct sums behave in the obvious way
so this is a complete description of the finite dimensional 
situation.

A physical setup realising these kinematics would have some
surprising properties.

First of all the map $\xi\otimes \eta \mapsto \eta\otimes\xi$
does not pass to the Connes tensor product so fermionic or bosonic
statistics would not make sense. It is known however that in
many examples of systems of bimodules there is a unitary
braiding which could be interpreted as playing the role of
an exchange of particles.

Secondly, the only observables for ${\cal H}_1$ that pass to 
the Connes tensor product are those which commute with the fusing
algebra $M$. In particular if the fusing algebra is non-commutative
the only fusing observables that remain observable on the combined system are those in
the centre of $M$.

Thirdly, if two sysems are so entangled that \underline{every} 
observable of one were equivalent to an observable for the
other then the Hilbert space for the fused system would collapse
to a one-dimensional one. This is reminiscent of the Pauli
exclusion principle.

\5
\subsection{Principal graphs.}
In the subfactor context it is natural to ask about the 
algebra generated by the $E_i$'s as above which actually
has a von Neumann algebra structure. The algebra generated by the first
$n$ $E_i$'s is finite dimensional and naturally 
included in the next one.  A very visual way to describe inclusions
of such finite dimensional algebras is by a ``Bratteli diagram'' (\cite{bratteli:af} which
records the ranks of the minimal projections of the smaller algebra
in the simple components of the bigger one. The tower of algebras of
the $E_i$'s as above has the Bratteli diagram  below (for index $\geq 4$).
\5
\label{ole}
\qquad \vpic{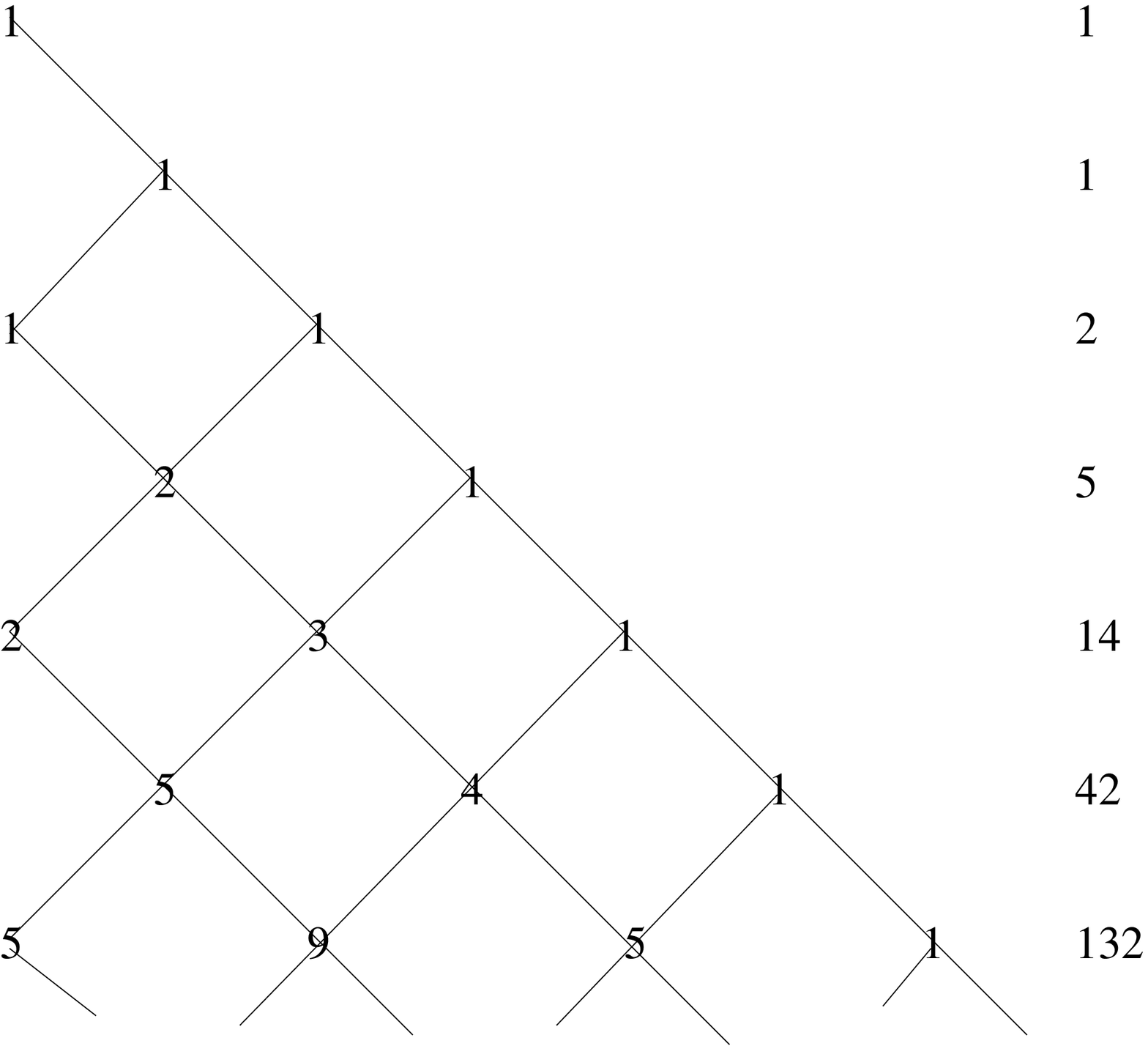} {4in}
\5
The numbers on the diagram are the sizes of the matrix algebras
which are the simple components,  and the numbers
to the right are the dimensions of the whole Temperley-Lieb algebra-
the Catalan numbers ${1\over n+1}{ 2n \choose n}$.

A provocative connection occurs here between this and section 
\ref{decomposition}.
For if we write $F_i=1+S_{i(i+1)}$ then these $F_i$ 
satisfy the relations \ref{tl1},\ref{tl2} and \ref{tl3} for $Q=4$.
So that the algebra generated by the $F_i$ should have (and indeed has
for index $\geq 4$) , the same Bratteli diagram as \ref{ole}!
Thus the decomposition of the tensor powers of the 2-dimensional representation
of $SU(2)$ are also described by \ref{ole}. We see that in fact
\ref{ole} is redundant in the sense that its essential information is
the graph $A_\infty$ of \ref{cgr}. The numbers on the Bratteli diagram are
just the number of loops on the graph $A_\infty$ starting and ending at 
the leftmost vertex. These loops form a basis of the algebra 
that gives the Bratteli diagram.

But there is a lot more finite dimensional algebra inside the tower $M_i$.
In fact the \underline{centralisers} $N'\cap M_i= \{x\in M_i |xn=nx \quad \forall n\in N\}$
are all finite dimensional and obviously each one is included in the
next. So they also have a Bratteli diagram which can be shown to have
the same structure as that for the $E_i's$ - there is a graph $\Gamma$ with
a privileged vertex *, such that the algebra $N'\cap M_i$ is given by
loops on $\Gamma$ based at *. The graph $\Gamma$ is called the 
\underline{principal graph}. There is a duality between $N\subseteq M$ and
$M\subseteq M_1$ and the principal graph of $M\subseteq M_1$ is called
the \underline{dual principal graph}.

There are subfactors for which the principal and dual principal graphs
are both $A_\infty$. More interestingly perhaps the subfactors in \cite{jones:index} which
give the index values $4\cos ^2 \pi/n$ have principal graphs $A_{n-1}$
(with one of the end points being * ). We will return to this 
when we discuss finite ``quantum'' subgroups of $SU(n)$.

For another example choose an outer action of a finite group $G$
on a II$_1$ factor $M$ and let $N=M^G$, the subfactor of fixed
points under the action of $G$. Then the principal and dual principal
graphs $\Gamma$ and $\check \Gamma$ of $M^G\subseteq M$ are as follows:
$\Gamma$ has a vertex * and a vertex for each irrep of $G$, with as
many edges between * and the irrep as the dimension of the irrep, and no
other edges. $\check \Gamma$ has a vertex * and one vertex for each element
of $G$ with an edge between * and each element of $G$, and no other edges.
Thus for the symmetric group $S_3$ the graphs are as below:
\5
\hspace{1in} \vpic{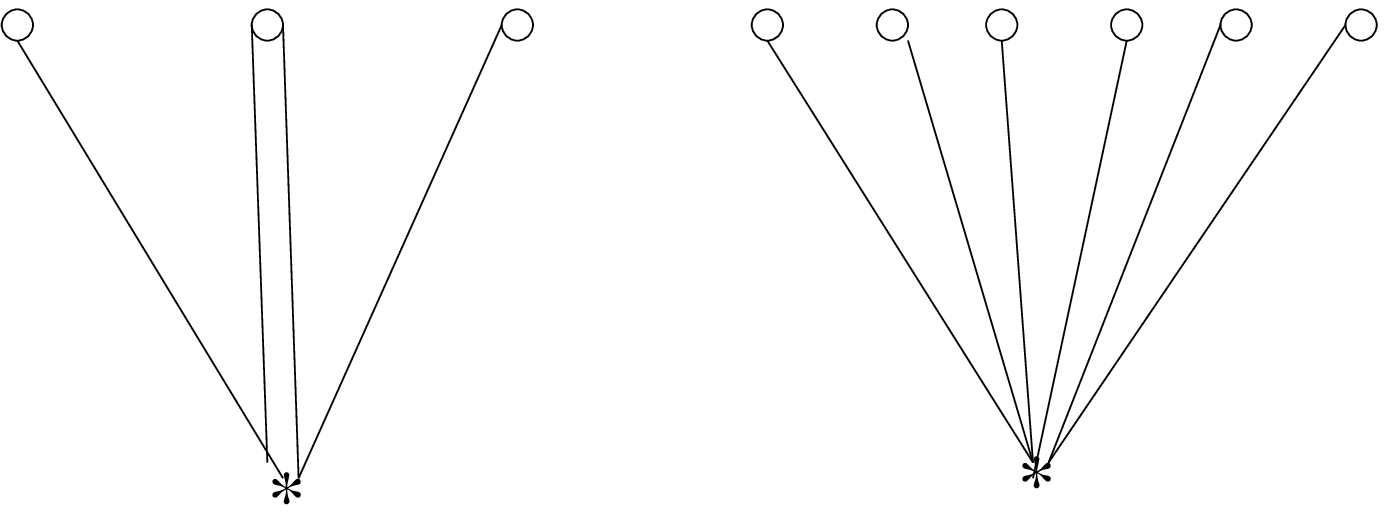} {3in}
\5

 \section{Braid group representations.}

 The braid group $B_n$ is the finitely presented group with presentation
$$<\sigma_1, \sigma_2,...,\sigma_{n-1}|\sigma_i\sigma_{i+1}\sigma_i =
\sigma_{i+1}\sigma_i\sigma_{i+1}\quad\rm{   for  }\quad i=1,2,..,n-2,$$
$$\hspace {0.7in}\rm{ and}\quad \sigma_i\sigma_j=
\sigma_j\sigma_i \quad for\quad |i-j|\geq 2>$$

If one puts $\lambda=\mu=\rho$ in the Yang Baxter equation \ref{ybe} it
is clear that one obtains a braid group representation 
by sending $\sigma_i$ to $R_{i(i+1)}(\lambda)$ on $\otimes^n V$ where $R(\lambda)$
acts on $V\otimes V$, provided $R(\lambda)$ is invertible. 
This was first done in \cite{jones:poly} for the matrix of \ref{rmatrix}. The resulting
braid group representation is of considerable interest. It may 
well be faithful for generic $q$. This method of obtaining braid group 
representations was applied universally after the development of
quantum groups. It was shown by Krammer and Bigelow (\cite{krammer},\cite{bigelow}) that
certain of the ensuing finite dimensional representations of
the braid group are indeed faithful.

The main reason for the interest in the braid group is that it
has geometric interpretations. First, according to its name, it
is the group of all braids on $n$ strings. A braid is a way 
of tying $n$ points on a bar at the top to $n$ points on 
a bar at the bottom, by strings whose tangent vector always has
a non-zero vertical component. Thus the figure below represents a
braid on 3 strings.
\5
\hspace{1.2in}\vpic{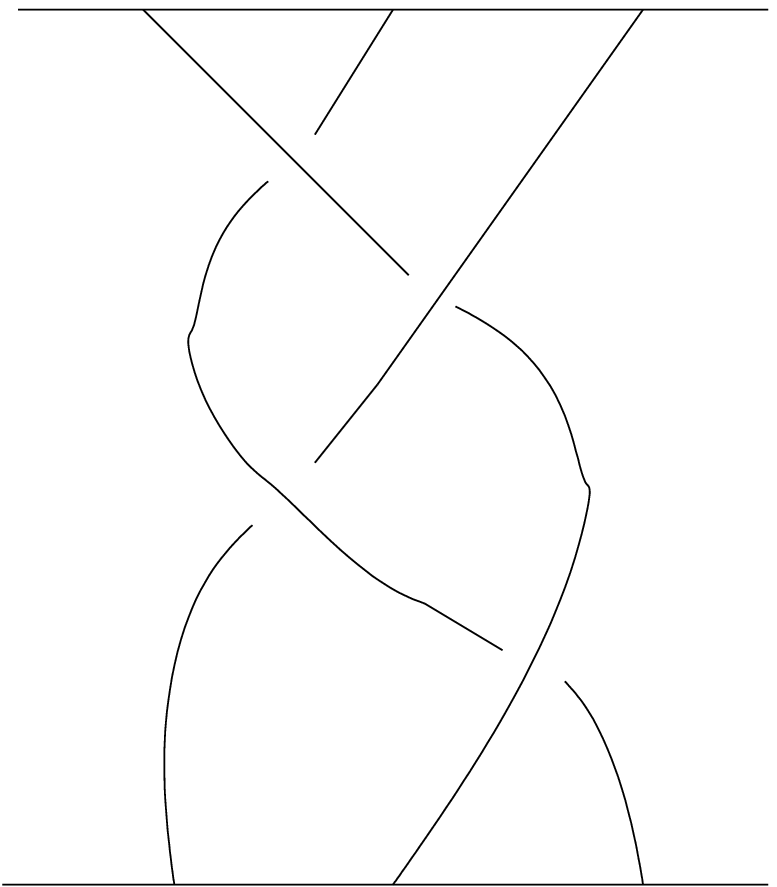} {1.5in} 
\5
Alternatively, braids can be thought of as motions of $n$ points
in the plane. The position at any time $t$ of the $n$ points
being determined by the intersection of a horizontal plane 
with $y$ co-ordinate equal to $t$ with the strings of the braid.
In this way $B_n=\pi_1((\mathbb C ^n\setminus \Delta) /S_n)$ where
$\Delta$ is the set of $n$-tuples $(z_1,z_2,...,z_n)$ with $z_i=z_j$
for some pair $i\neq k$ and the symmetric group $S_n$ acts in the
obvious way.

Knots and links can be formed from braids by tying the tops of the
strings to the bottom. The figure-eight knot is obtained by doing
this to the braid drawn above. It was realised in \cite{jones:poly} that the
trace coming from the subfactor origin of the braid group 
representation furnished an invariant of knots and links.

From the physics point of view the most interesting property
of these representations is their unitarity. This is a difficult 
topic in general for the representations are not unitary for
real positive values of $q$, even though there is a natural
Hilbert space structure on the vector spaces on which they act.
Fortunately this Hilbert space structure persists
enough to supply, for fixed $n$, a small interval of $q$ values 
(containing $1$) on the circle
for which the representation is unitary on the Hilbert space
$\otimes^n V$ though the Hilbert space structure will fail if
$n$ is increased indefinitely and $q$ is left fixed.
For $q$ a root of unity of the form $e^{2\pi i \over p}$ 
there are other statistical mechanical models with Hilbert
spaces for which the braid group representations are unitary
for all $n$ and fixed $q$. These are also the values of $q$
for which subfactors of finite depth occur - see Wenzl - \cite{wenzl} and Xu -\cite{feng}.
Subfactors can be constructed for positive real values of 
$q$ as well. They are of infinite depth  and  have been analyzed by Sawin in \cite{sawin}.

\section{Detective work}
Hopefully the reader has been struck by common threads, notational and
otherwise, connecting all the previous sections. They suggest a grand 
structure in which the formulae we have come up with are a small but
significant part. We are not convinced that the last word has been said
on this grand structure but for the case of the vertex models the notion
of quantum group does the trick with great elegance and power.

Beginning with the decomposition of the tensor powers of the representations
of $SU(2)$ we have presented operators that ``deform'' the permutation
operators $S_{i(i+1)}$ of \ref{decomposition}, the most general of which
was the $R$-matrix \ref{rmatrix}. We have only hinted that the theory
goes beyond $SU(2)$ but in fact Cherednik gave $R$-matrices that
deform the representation of the symmetric group on the tensor powers
of $\mathbb C^k$ for $k\geq 2$ and Wenzl independently discovered 
the same objects in the subfactor/braid group context in \cite{wenzl2}. Thus it
was natural to hope for an object that would ''deform'' $SU(k)$ in
a way that Schur-Weyl duality would be preserved and the commutant of
this object would be generated by the $R$-matrices.
This was done by Jimbo and Drinfeld who were greatly inspired also by
a vision of this procedure as a ``quantisation'' of the theory of
integrable systems in Hamiltonian mechanics.  See \cite{drinfeld},\cite{jimbo} and \cite{ft}.

There are many accounts of this work and we do not want to dwell on
it as we are really interested here in cases where the quantum 
group formalism does not apply easily but which arise naturally in
the subfactor world. Suffice it to say that the final result is
the construction of $R$-matrix solutions (depending in various 
ways on the spectral parameter) for all (finite dimensional)simple Lie algebras and
all representations thereof. In the statistical mechanics formalism
the horizontal and vertical directions on the lattice may 
correspond to different representations of the Lie algebra. And there
are extensions to affine Lie algebras.

Appropriate $R$-matrices can be evaluated at special values
of the parameters so that they give braid group representations
and all such representations are known to give link invariants
where a representation of the Lie algebra can be assigned independently
to each component of the link. The invariants are always polynomials
in the quantum deformation parameter $q$.
See \cite{rosso}. The link invariants are powerful but many elementary 
questions remain unanswered. Perhaps the most galling of these is
the question of whether the simplest of the invariants,
corresponding to the Lie algebra $sl(2)$ and its $2$-dimensional
representation, detects knottedness, i.e. is there a non-trivial knot for 
which this polynomial is the same as for the unknot?
The answer is known for links (see \cite{ekt}) and all knots up to 
17 crossings have been checked.

Kohno showed in \cite{kohno} that the braid group representations coming
from quantum groups could be obtained using only the data of
the ``classical'' Yang-Baxter equation (as formulated by Drinfeld)
which could be used to define a flat connection on configuration
space $(\mathbb C ^n\setminus \Delta) /S_n$.

\section{The Haagerup and Haagerup-Asaeda subfactors.}

In the light of the huge success of quantum groups and connections
with conformal field theory and $\sigma$-models which we have not
mentioned here, one is entitled to ask if all interesting theories
are somehow obtainable from Lie groups and ``geometry''. This 
is an important question, particularly for subfactors whose 
utility would be brought into question if one could obtain all
examples from structures outside subfactor theory. In this section
and the next we offer two examples of how subfactors can produce
examples independently of any geometric input, the first pioneered
by Haagerup and the second by Ocneanu.

Haagerup asked the question: what is the (irreducible) subfactor of smallest
index greater than $4$ that occurs for a II$_1$ factor? In fact
Popa has shown that any index value greater than $4$ occurs but those
examples have $A_\infty$ as principal graph. So Haagerup's real question
was to find the smallest index subfactor with principal graph different
from $A_\infty$. This he solved in spectacular fashion. He first showed
in \cite{haagerup} that the smallest index value is $5+\sqrt{13}\over 2$ and identified
the principal graph and dual principal graph as being one of the two 
below (which are dual to each other):
\5
\vpic{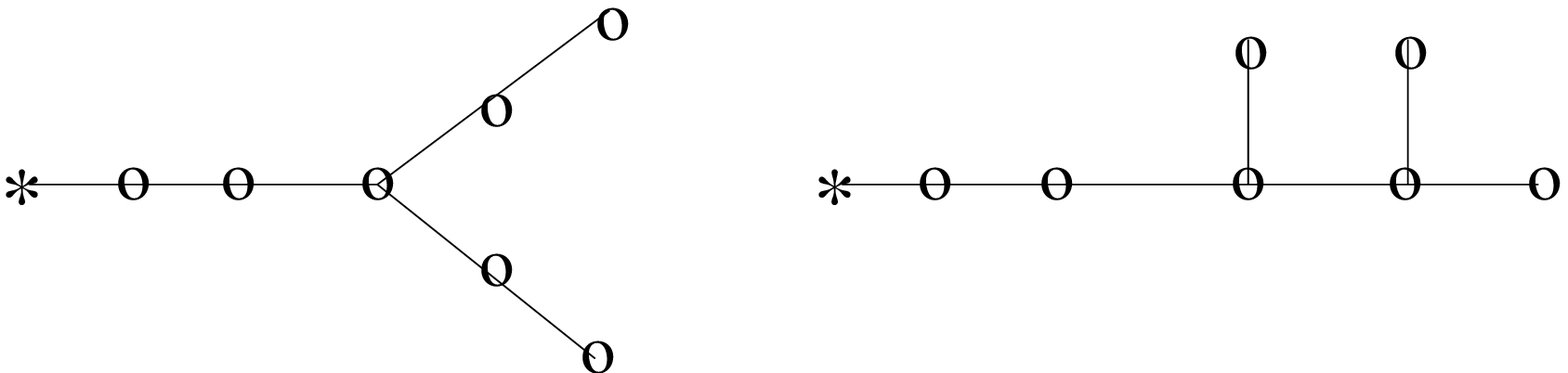} {3.5in}
\5

Then in \cite{haagas} Haagerup and Asaeda showed that there is indeed a subfactor
of the hyperfinite II$_1$ factor which has these graphs as dual principal
graphs. They further showed that the next possible index value is  $5+\sqrt{17}\over 2$,
constructed the subfactor and calculated the principal graphs.

The techniques for eliminating other graphs and constructing
the examples were highly calculatory, relying on Ocneanu's
theory of connections. It is a major challenge in the theory
to come up with an interpretation of these subfactors  
subfactors as members of a family related to some other mathematical
objects.  Izumi in \cite{izumi} made some good progress in this direction.
In \cite{jones:quadratic} we have taken a diagrammatic approach and calculated certain
parameters that show that the subfactors of indices  $5+\sqrt{13}\over 2$
and  $5+\sqrt{17}\over 2$ are of a very different kind.
\section{Finite ``quantum'' subgroups of Lie groups.}

 Subfactors of index less than 4 have principal graphs equal to
one of the A,D or E Coxeter-Dynkin diagrams. (Ocneanu showed that in
fact only the diagrams $D_{2n}$ and $E_6$ and $E_8$ can actually occur-
and he constructed subfactors for each of these diagrams.)
We saw in the first section how the extended Coxeter-Dynkin diagrams
are connected with finite subgroups of $SU(2)$. There is also a 
subfactor connected with a finite subgroup of $SU(2)$ as follows:
Construct a II$_1$ factor $R$ as $\otimes^\infty M_2(\mathbb C)$
(completed using the normalised trace to get a von Neumann algebra). 
The subfactor $R_0$ is the subalgebra of all elements of $R$ of
the form $id\otimes x$ where $id$ is the $2\times 2$ identity matrix.
The group $SU(2)$ acts on $R$ by the infinite tensor product of
its action on $M_2(\mathbb C)$ by conjugation. This action preserves
the subfactor $R_0$. So for any subgroup $G$ of $SU(2)$ one may
form the subfactor $N\subseteq M$ of fixed points for $G$: $N=R_0^G$
and $M=R^G$. The principal graph for $N\subseteq M$ is then
the extended Coxeter-Dynkin graph for the subgroup! The index is $[M:N]=4$.

This suggests that one should be able to interpret the vertices of the
principal graph as representations of something and the edges
as induction/restriction. This is indeed possible and is inevitable
in Connes' picture where bimodules over a II$_1$ factor replace
representations of a group. The vertices of the principal graph
are a certain family of $N-N$ bimodules and $N-M$ bimodules and
the edges count induction/restriction multiplicities. This point
of view was first pointed out  by Ocneanu.

Putting together the index $<4$ and index $4$ cases above we
see that it is natural to think of the index $<4$ subfactors as
being some kind of quantum version of subgroups of $SU(2)$
whose ``representation theory'' is a truncation of the representation
theory of the corresponding genuine subgroup of $SU(2)$.

Cappelli, Itzykson and Zuber in \cite{ciz} also ran into the $A-D-E$ Coxeter-Dynkin
diagrams in their attempt to understand modular invariants in 
conformal field theory. In extending that work DiFrancesco and Zuber in \cite{dfz}
looked for truncations of the representation theory graphs of 
subgroups of $SU(3)$ that could give modular invariants. Zuber 
presented a list of such graphs satisfying certain criteria and
conjectured that it was complete. Ocneanu used the subfactor
point of view to exhibit the complete list for $SU(3)$ and beyond,
with slight changes in the requirements for a graph to be on the 
list. Unfortunately the situation is a little too nice for $SU(2)$ 
because all its representations are equivalent to their conjugates.
To do justice to the subfactor point of view would require a detour
beyond the principal graphs so we simply refer to Ocneanu's paper \cite{ocneanu:bariloche}
for those interested in this topic.

\section{The direct relevance of subfactors to physics.}
The connection between subfactors and physics that we have outlined
above has been somewhat indirect, passing from certain elements in
centraliser towers to quantum spin chains and/or statistical mechanical
models. However probably von Neumann's main motivation for studying
his algebras was because of their relevance to the quantum mechanical 
formalism. So one might hope for a more direct connection between
subfactors and quantum physics. This does exist and is associated with
many names. It goes back to the pioneering work of Haag and Kastler
who sought to develop a non-perturbative framework for quantum field
theory by taking as basic ingredients the algebras of observables
localised in various regions of space-time.

I feel unable to give a satisfactory account of this theory and it
is to be hoped that a book on it will appear in the near future.
I will just say the following-a subfactor naturally appears by
looking at the von Neumann algebras associated to certain regions
of space-time. If two regions of space time are such that no
light ray can connect them (they are not causally connected),
 then their von Neumann algebras of local
observables should commute. These von Neumann algebras are known to
be factors (of type III) so one may obtain a subfactor by taking
$N$ to be the algebra of observables localised in a region and
$M$ to be the commutant of all observables localised in the causal
complement of that region.

Using this and other motivations, Wassermann and I looked for 
subfactors in the theory of loop group representations. A unitary
(projective) representation of the group $LSU(2)$ of smooth functions from
the circle to $SU(2)$ can be thought of as a one-dimensional quantum field theory
whose currents are given by functions from the circle into the Lie algebra $su(2)$
and whose algebra of observables localised in an interval $I$ of the
circle is the von Neumann algebra generated by the normal subgroup
$L_ISU(2)$ of loops supported in
that interval. In this way we were led to the subfactor (where $I^c$ is 
the interval complementary to $I$ on the circle) 
$$L_ISU(2)''\subseteq L_{I^c}SU(2)'.$$
Wassermann subsequently showed in \cite{wassermann} that the set of indices of subfactors so obtained does 
indeed contain the set $\{4\cos^2\pi/n\}$ and extended this work
to $SU(n)$, the diffeomorphism group of the circle,  and beyond.

\bibliographystyle{amsplain} \providecommand{\bysame}
{\leavevmode\hbox to3em{\hrulefill}\thinspace}

\end{document}